%% file: final-for-arxiv.tex
\documentclass[a4paper]{amsart}
\usepackage{amsthm,amssymb}
\usepackage[T1]{fontenc}
\usepackage{enumerate}
\usepackage{mathrsfs}
\usepackage{dsfont}
\usepackage{mathtools}
\usepackage{scalerel,stackengine}
\usepackage{nicefrac}
\usepackage{natbib}

\usepackage{caption}
\usepackage[colorlinks,citecolor=blue,urlcolor=blue]{hyperref}

\usepackage{tikz,tikz-cd}
\usetikzlibrary{calc,hobby}
\usetikzlibrary{shapes}
\tikzset{point/.style={fill=black,circle,inner sep=1pt}}

\numberwithin{equation}{section}

\newcommand\conr{\con_R} %the frame contact relation
\newcommand\conp{\con_P} %the partition contact relation
\newcommand\conm{\con_{\diamond}} %contact from a modal operator
\newcommand\nconm{\notcon_{\diamond}} %contact from a modal operator - negation
\newcommand{\Cfourc}{\boldsymbol{\mathsf{C4^c}}} %the category of complete BCAs with (C4c)

\input{basic_mac.tex}

\renewcommand{\Nat}{\omega}

\input{bbm-without-txfonts.tex}

\title{From contact relations to modal operators, and back}

\author[R. Gruszczy\'nski and P. Mench\'on]{Rafa\l\ Gruszczy\'nski and Paula Mench\'on}\thanks{This research was funded by (a) the National Science Center (Poland), grant number~2020/39/B/HS1/00216 and (b) the European Union’s Horizon 2020
research and innovation program under the Marie Skłodowska-Curie grant agreement
No 101007627. For the purpose of Open Access, the authors have applied a CC-BY public copyright license to any Author Accepted Manuscript (AAM) version arising from this submission.}

\date{}

\address{Rafa\l\ Gruszczy\'nski, Paula Mench\'on\\
Department of Logic\\
Nicolaus Copernicus University in Toru\'n\\
Poland}

\email{gruszka@umk.pl, paula.menchon@v.umk.pl}
\begin{document}

\begin{abstract}
  One of the standard axioms for Boolean Contact Algebras says that if a region $x$ is in contact with the join of $y$ and $z$, then $x$ is in contact with at least one of the two regions. Our intention is to examine a~stronger version of this axiom according to which if $x$ is in contact with the supremum of some family~$S$ of regions, then there is a~$y$ in $S$ that is in contact with~$x$. We study a~modal possibility operator which is definable in complete algebras in the presence of the aforementioned axiom, and we prove that the class of complete algebras satisfying the axiom is closely related to the class of modal KTB-algebras. We also demonstrate that in the class of complete extensional contact algebras the axiom is equivalent to the statement: \emph{every region is isolated}.  Finally, we present an~interpretation of the modal operator in the class of the so-called \emph{resolution contact algebras}.
  
  \smallskip
  
  
\end{abstract}

\maketitle

\section{Introduction}

Since the beginning of the twentieth century Boolean Contact Algebras have been the standard algebraic approach to region-based analyses of spatial relations of contact (proximity) and separation. The metatheory of the algebras have been developed, and various subclasses and extensions of the class of the algebras have been studied. In this paper we intend to contribute to the field via an analysis of the stronger form of one of the standard axioms for Boolean Contact Algebras. On the intuitive level, the standard axiom says that if a~region $x$ is in contact with the Boolean meet of regions $y$ and $z$, then $x$ is in contact with either $y$ or $z$ (or, possibly, both). The strengthening we have in mind generalizes the axiom to infinite cases; that is, the axiom we put forward says that if $x$ is in contact with the supremum of the family of regions $S$, then in $S$ there is a~region that is in contact with~$x$. 

The axiom is worth studying due to at least two reasons. Firstly, \cite{Duntsch-et-al-RBTODSAPA} proved that every Boolean contact algebra can be embedded in a particular type of algebra that satisfies our axiom, the so-called relational algebra. Secondly, in the class of complete Boolean Contact Algebras that satisfy the axiom, we can define a modal possibility operator. This establishes a close-knit connection between modal algebras and a subclass of Boolean Contact Algebras.

% An interesting property is that due to a result by \cite{Duntsch-et-al-RBTODSAPA}, every contact algebra can be embedded in a~ relational algebra, a complete Boolean algebra with a contact relation that satisfies our stronger axiom. 

In Section~\ref{sec:preliminaries} we recall some standard notions and their properties necessary for our results. The notions include \emph{grills} of Boolean algebras, \emph{modal operators} that can be obtained from grills and ideals, and a~\emph{quasi-modal operator} of Celani's (\citeyear{Celani-QMA}). We define the notion of a~\emph{co-principal grill} and by means of it we define a~possibility modal operator.

Section~\ref{sec:BCAs} recalls elementary facts about Boolean Contact Algebras, introduces a quasi-modal operator related directly to the contact relation, and summarizes properties of the subordination relation, which plays an important role in the remaining sections.

In Section~\ref{sec:examples} we provide various examples and counterexamples whose aim is to familiarize the reader with the content of the new axiom. Among others, we show that the axiom is indeed stronger than the standard one, that every relational contact algebra satisfies it, and that under the standard topological interpretation of contact the regular closed algebra of any Alexandroff space satisfies it too.

Section~\ref{sec:extensional} contains a proof of the fact that in the so-called \emph{extensional} Boolean Contact Algebras introduction of the axiom results in collapsing the contact relation to the standard overlap relation. This is a situation equivalent to every region being isolated, i.e., every region being its own non-tangential part. 

In section~\ref{sec:the-operator} we study mutual dependencies between the contact relation and a~modal possibility operator. Firstly, we demonstrate that if a Boolean contact algebra satisfies the axiom, then we can define a~modal operator that is completely additive and we investigate its properties. Secondly, we show that in ever complete KTB-algebra we can define a~contact relation that satisfies the axiom.
Thirdly, we gather Boolean Contact Algebras satisfying the axiom into a~category using the so-called \emph{p-morphisms}, and we prove that this category is isomorphic to the category of the complete KTB-algebras and homomorphisms.

The last section introduces the class of resolution algebras which serve as a~kind of geometric interpretation of Boolean Contact Algebras satisfying the axiom. We show that there is a~dependence between these algebras and the class of S5 modal algebras; that is, an expansion of any such an algebra with the contact relation can be always embedded into a~modal expansion of a~resolution algebra.

\section{Preliminaries}\label{sec:preliminaries}
Let $\frB=\langle B,\mathord{\cdot},\mathord{+},-,\zero,\one\rangle$
be a~Boolean algebra (BA for short) with the operations of, respectively, meet, join and boolean complement; and with the two distinguished elements: the minimum $\zero$ and the maximum $\one$. Elements of the domain are called \emph{regions}. We use the standard logical operators and connectives: $\neg$ (negation), $\wedge$ (conjunction), $\vee$ (disjunction), $\rarrow$ (material implication), $\iff$ (material equivalence), $\exists$ (existential quantifier) and $\forall$ (universal quantifier). $\iffdef$ is a definitional equivalence, and $\defeq$ a~definitional equality. If $S\subseteq B$, then $\bigvee S$ is the supremum of $S$ (its least upper bound) and $\bigwedge S$ is its infimum (greatest lower bound). Given a~domain $D$ and its subset $A$, $A^\complement$ is the set-theoretic complement of~$A$ with respect to~$D$. Throughout the paper, `iff' is an abbreviation for `if and only if'.

% $\bigvee_{i\in I} x_i$ is an operation of the indexed supremum, where $I$ is a (possibly infinite) set of indexes. Thus, $\bigvee_{i\in I} x_i$ is the least upper bound of the set $\{x_i\mid i\in I\}$.\marginpar{\tiny Are indexed families necessary now after Paula's changes? Check this.} The analogous convention is used for indexed infima. 

%We adopt the convention according to which we denote algebras by means of fraktur letters, and their domains by means of their Latin counterparts.

The class of all BAs is denoted by `$\BA$', and the class of all \emph{complete} BAs by `$\BAc$'. %Boolean lagebras will be denoted by fraktur letters $\frA$, $\frB$ and so on, and we introduce the convention that their Latin counterparts will serve as names of the domains of suitable algebras. 

As usual, in $\frB\in\BA$ we define two standard order relations:
\begin{align*}
x\leq y&{}\iffdef x\cdot y=x\,,\\%\tag{$\mathrm{df}\,\mathord{\leq}$}\\
x<y&{}\iffdef x\leq y\wedge x\neq y\,.%\tag{$\mathrm{df}\,\mathord{<}$}
\end{align*}
%In the former case we say that $x$ is \emph{part} of~$y$ or that $x$ is \emph{below}~$y$, in the latter that $x$ is a~\emph{proper part} of~$y$ or that $x$ is \emph{strictly below}~$y$. 
We call $\scrS\subseteq B$ an \textit{upward closed subset} if it satisfies the following condition for all $x\in B$:
\[
x\in\scrS\rarrow\upop x\subseteq\scrS\,,
 \]
 where:
 \begin{equation*}%\tag{$\dftt{\uparrow}$}
     \upop x\defeq\{y\in B\mid x\leq y\}\,.
 \end{equation*}
Analogously, $\scrS\subseteq B$ is a \textit{downward closed subset} if it satisfies the following condition for all $x\in B$:
\[
x\in\scrS\rarrow\downop x\subseteq\scrS\,,
 \]
 where:
 \begin{equation*}%\tag{$\dftt{\uparrow}$}
     \downop x\defeq\{y\in B\mid y\leq x\}\,.
 \end{equation*}

Additionaly, we also define a~binary relation of \emph{overlapping} or \emph{compatibility} of regions:
\begin{equation*}%\tag{$\dftt{\Overl}$}
    x\Overl y\iffdef x\cdot y\neq\zero\,,
\end{equation*}
and its complement, the  \emph{disjointness} (or \emph{incompatibility}) of regions: 
\begin{equation*}%\tag{$\dftt{\ext}$}
    x\ext y\iffdef x\cdot y=\zero\,.
\end{equation*}
For any region $x$ we define:
\begin{equation*}
    \Overl(x)\defeq\{y\in B\mid y\Overl x\}\,,
\end{equation*}
the set of all regions that overlap $x$. Every Boolean algebra is $\Overl$-extensional, i.e., it satisfies the following condition for all $x,y\in B$:
\[
\Overl(x)=\Overl(y)\rarrow x=y\,.
\]

As usual, a subset $\fil \subseteq B$ is a \emph{filter}
of $\frB$ if it is a non-empty upward closed set such that for all $a, b \in \fil$, $a \cdot  b \in \fil$. Analogously, a subset $\scrJ \subseteq B$ is an \emph{ideal}
of $\frB$ if it is a non-empty downward closed set such that for all $a, b \in \scrJ$, $a +  b \in \scrJ$.
\subsection{Completely prime grills} We recollect some elementary properties of Boolean Contact Algebras. For more details see \citep{Dimov-et-al-CARBTSPA1,Dimov-et-al-CARBTSPA2,Vakarelov-et-al-PATSRBTS,Vakarelov-et-al-ANOPSACBM}, and \citep{Naimpally-et-al-PS}. 

% By a \emph{stack} of a~Boolean algebra we understand any non-empty subset $\scrS$ of its domain that does not contain the zero element and is upward closed:
% \[
% x\in\scrS\rarrow\upop x\subseteq\scrS\,,
% \]
% where:
% \begin{equation}\tag{$\dftt{\uparrow}$}
%     \upop x\defeq\{y\in B\mid x\leq y\}\,.
% \end{equation}
A non-empty subset $\scrS$ of the domain of a~Boolean algebra is a \emph{grill} iff it does not contain the zero element, it is upward closed and it also satisfies the following condition for all $x,y\in B$:
\[
x+y\in\scrS\rarrow x\in\scrS\vee y\in\scrS\,.
\]
It is routinely verified that the set-theoretic complement of a grill is an ideal. We denote by $\Gr(\frB)$ the set of all grills of a Boolean algebra~$\frB$. A grill  $\scrG$ is \emph{completely prime} if it satisfies the following, stronger property:
\[
\textstyle\bigvee J\in\scrG\rarrow(\exists x\in J)\,x\in\scrG
\]
for all arbitrary families of regions $J\subseteq B$ such that their join exists.
For example, for any non-zero region $x$, $\Overl(x)$ is a completely prime grill.
\begin{proposition}\label{prop:cpGrills-principal-ideals}
Let $\frB$ be a complete Boolean algebra. Then, a grill $\scrG$ is completely prime iff $\scrJ\defeq \scrG^\complement$ is a~principal ideal, i.e., $\scrJ=\downop x$ for some $x\in B$.%, where\/\textup{:}
%\begin{equation*}%\tag{$\dftt{\downop}$}
%\downop x\defeq\{y\in B\mid y\leq x\}\,.
%\end{equation*}
\end{proposition}
\begin{proof}
If $\scrG$ is a completely prime (c.p., in abbreviation) grill, then $\scrJ\defeq\scrG^\complement$ is an ideal that satisfies the following completeness property:
\[
\textstyle J\subseteq\scrJ\rarrow\bigvee J\in\scrJ\,.
\]
But then $\scrJ=\downop\bigvee\scrJ$. The other way round, suppose that $\scrS=\downop x$  for some $x\in B$ and let $J\subseteq B$ be a family of regions such that $\bigvee J \notin\downop x$. Then, there must be at least one $y\in J$ such that $y\nleq x$.
\end{proof}
As a consequence, given a complete Boolean algebra, $\scrG$ is a c.p. grill iff there exists a~region $x$ such that $\scrG=(\downop x)^\complement=\Overl(-x)$. 

% We denote by $\cpGr(\frB)$ the set of all c.p. grills of a complete Boolean algebra~$\frB$:
% \begin{equation*}%\tag{$\dftt{\cpGr}$}
%     \cpGr(\frB)\defeq\{\Overl(x)\mid x\neq\zero\}\,.
% \end{equation*}
The above results show that the notion of a c.p. grill is non-trivial, in the sense that c.p. grills exist in every non-degenerate algebra (i.e., an algebra with at least two elements).

\begin{corollary}
  There is a~one-to-one correspondence between c.p. grills and non-zero regions of a~given complete BA, and so between c.p. grills and proper principal ideals.
\end{corollary}

It is well known that every grill is equal to the set-theoretic sum of all ultrafilters that are its subsets, and it is easy to verify that any sum of ultrafilters is a~grill \citep{Thron-PSAG}. Therefore, it might be tempting to identify completely prime grills with the sums of completely prime filters.\footnote{Recall that a~filter $\fil$ is called \emph{completely prime} if for all families $J\subseteq B$ such that $\bigvee J\in\fil$ there exists $x\in J$ satisfying $x\in\fil$.} However, this cannot hold in general in light of the following  well-known property of BAs:
\begin{proposition}[\citealt{Picado-Pultr-FL}]
In every complete BA, $\fil$ is a c.p. filter iff there exists an atom\footnote{A~region $x$ of a Boolean algebra $B$ is an \emph{atom} iff $\zero<x$ and $x$ is minimal with respect to~$<$.} $a$ such that $\fil=\upop a$ (i.e., c.p. filters are exactly principal ultrafilters\footnote{An ultrafilter~$\ult$ of a BA is \emph{principal} iff $\ult={}\upop a$, for some atom~$a$.}).
\end{proposition}
\begin{proof}
Indeed, if $\fil$ is a c.p. filter, then $\bigwedge\fil\in\fil$. Since if not, $\bigvee\{-x\mid x\in\fil\}\in\fil$, so there exists $x\in\fil$ such that $-x\in\fil$, a contradiction. Thus $\fil=\upop\bigwedge\fil$, and since $\fil$ is an ultrafilter, then $\bigwedge\fil$ must be an atom.
\end{proof}
Therefore, atomless algebras do not contain completely prime filters, yet they contain c.p. grills, as we can always construct principal ideals. What we can only assert is that:
\begin{proposition}
If a BA contains c.p. filters, then any sum of these is a~completely prime grill.
\end{proposition}

The situation is different in the case of complete and atomic algebras. If $x\in\scrG$, where $\scrG$ is a~c.p. grill, then $x\neq\zero$, so there is a~set of atoms $S\neq\emptyset$ such that $x=\bigvee S$. In consequence there is an atom $a\in S\cap\scrG$, and $\upop a\subseteq\scrG$. Thus we have:
\begin{proposition}
If $\frB\in\BAc$ is atomic, then for any c.p. grill $\scrG$\textup{:}
\begin{enumerate}
    \item for every $x\in\scrG$, there is a~principal ultrafilter $\ult\subseteq\scrG$ such that $x\in\ult$,
    \item $\scrG$ contains at least one atom,
    \item $\scrG$ is the sum of all principal ultrafilters contained in $\scrG$.
\end{enumerate}
\end{proposition}

% \begin{definition}
% Given a family $\mathbf{F}$ of sets of regions of a Boolean algebra, we will say that $\mathbf{F}$ is \emph{saturated} iff for every non-zero region $x$ there is an $F\in\mathbf{F}$ such that $x\in F$.
% \end{definition}

% \begin{corollary}
% The family of c.p. grills of any complete BA is saturated.
% \end{corollary}
% \begin{proof}
% Obvious.
% \end{proof}

\subsection{Modal operators from ideals and grills}

Let us first recall some standard notions.
\begin{definition}
A \emph{possibility modal operator} on a BA $\frB$ is a function $\Diamond\colon B\fun B$ that is normal: $\Diamond(\zero)=\zero$, and additive: $\Diamond(x+y)=\Diamond(x)+\Diamond(y)$ for all $x,y\in B$. $\Diamond$ is \emph{completely additive} iff $\Diamond\left(\bigvee J\right)=\bigvee_{x\in J}\Diamond(x)$ for all families $J\subseteq B$ such that their join exists.

A \emph{necessity modal operator} is a function $\Box\colon B\fun B$ that is co-normal: $\Box(\one)=\one$, and multiplicative: $\Box(x\cdot y)=\Box(x)\cdot\Box(y)$ for all $x,y\in B$. $\Box$ is \emph{completely multiplicative} iff $\Box\left(\bigwedge J\right)=\bigwedge_{x\in J}\Box(x)$ for all families $J\subseteq B$ such that their meet exists.
\end{definition}

\begin{definition}
If $\Ide(\frB)$ is the lattice of ideals of a Boolean algebra $\frB$, then 
a~mapping $\Delta\colon B\to\Ide(\frB)$ is a~\emph{quasi-modal operator}\footnote{The notion was introduced by \cite{Celani-QMA}. } if it satisfies the following two conditions:
\begin{enumerate}
    \item $\Delta(\one)=L$,
    \item $\Delta(x\cdot y)=\Delta(x)\cap\Delta(y)$ for all $x,y\in B$.
\end{enumerate}
$\Delta$ is \emph{principal} iff $\Delta(x)$ is a~principal ideal for every $x\in B$, that is, for every $x\in B$ there exists $y\in B$ such that $\Delta(x)=\downop y$. If $\Delta$ is a principal operator, then---after \cite{Celani-SASIBDLWUO}---we can recover a~necessity operator $\Box\colon B\to B$ via:
\[
\Box x\defeq (\iota y)\, \Delta (x)=\downop y\,.\footnotemark
\]
\footnotetext{$\iota$ is the standard definite description operator, so $(\iota x)\,\varphi(x)$ is the only object $x$ that satisfies $\varphi$.}
Conversely, in a modal algebra $\langle \frB,\mathord{\Box}\rangle$ we can define a~quasi-modal operator $\Delta\colon B\to \Ide(\frB)$ by $\Delta(x)\defeq\downop\Box x$.
\end{definition}

Since grills are set-theoretic complements of proper ideals we can use them as well to define quasi-modal and modal operators. To this end we define:
\[
\Gr(\frB)^\circ\defeq\left\{\scrJ^\complement\suchthat \scrJ\in\Ide(\frB)\right\}=\Gr(\frB)\cup\{\emptyset\}\,,
\]
which is the lattice dual to $\Ide(\frB)$. Define an operator $\nabla\colon B\to\Gr(\frB)^\circ$ such that:
\[
\nabla(x)\defeq\Delta(-x)^\complement\,.
\]
It is routine to verify that:
\begin{equation}\label{eq:nabla-properties}
\nabla(\zero)=\emptyset\qquad\text{and}\qquad\nabla(x+y)=\nabla(x)+\nabla(y)\text{ for all }x,y\in B\,.
\end{equation}
So $\nabla$ is a quasi-modal possibility operator.\footnote{Observe that this a different operator than $\nabla$ from~\citep{Celani-QMA}.}

Let us say that a~grill $\scrG$ is \emph{co-principal} iff there is $x\in B$ such that $\scrG=(\downop x)^\complement$. Accordingly, $\nabla$ is \emph{co-principal} iff for every $x\in B$, $\nabla(x)$ is co-principal. As in the case of a~principal $\Delta$, for a~co-principal $\nabla$ we may define a~possibility operator $\Diamond\colon B\to B$ via:
\begin{equation}\label{df:Diamond}
\begin{split}
\Diamond x\defeq&{} (\iota y)\,\nabla(x)=(\downop -y)^\complement\\
=&{}(\iota y)\,\nabla(x)=\Overl(y)\,.
\end{split}
\end{equation}
%Since every grill is co-principal iff every ideal is principal, 
We can see that the operators $\Box$ and $\Diamond$ introduced by means of $\Delta$ and $\nabla$, respectively, are related to each other in the expected way for $x,y\in B$:
\begin{align*}
 \Diamond x=y&{}\iff \nabla(x)=(\downop-y)^\complement\\
 &{}\iff\Delta(-x)^\complement=(\downop-y)^\complement\\
 &{}\iff\Delta(-x)=\downop-y\\
 &{}\iff\Box-x=-y\\
 &{}\iff-\Box-x=y\,.   
\end{align*}
%As it can be seen, 
The above equivalence could be taken as the main definition based on which we can prove that \eqref{df:Diamond} holds.

We also have:
\begin{proposition}
If $\langle \frB,\Diamond\rangle$ is a~modal algebra, then $\nabla\colon B\to\Gr(\frB)^\circ$ such that $\nabla(x)\defeq(\downop -\Diamond x)^\complement$ is a quasi-modal operator.
\end{proposition}

The approach is entirely dual, in the following sense. Given a Boolean algebra~$\frB$, let us start with a quasi-modal operator $\nabla\colon B\to\Gr(\frB)^\circ$ that satisfies the properties from \eqref{eq:nabla-properties}. From $\nabla$ we can recover 
$\Delta$ via:
\begin{equation}
    \Delta(x)\defeq\nabla(-x)^\complement\,.
\end{equation}
Further, we can define a possibility operator $\Diamond$ as in \eqref{df:Diamond}, for a~$\Delta$ operator that is co-principal, and we can define a necessity operator. For such defined operators we can prove that all the suitable equivalences hold.

As the notion of \emph{grill} is closely related to Boolean Contact Algebras, whose subclass is the topic of this paper, in the sequel we will work with grills and possibilities rather than ideals and necessities.

\section{Boolean Contact Algebras}\label{sec:BCAs}

Any Boolean algebra is turned into a~\emph{Boolean contact algebra} by expanding it to a~structure $\langle\frB,\mathord{\con}\rangle=\langle B,\mathord{\cdot},\mathord{+},-,\zero,\one,\con\rangle$ where $\mathord{\con}\subseteq B^2$ is a~\emph{contact} relation which satisfies the following five universal axioms:
\begin{gather}
\neg(\zero\con x),\label{C0}\tag{C0}\\
x\leq y\wedge x\neq\zero\rarrow x\mathrel{\mathsf{C}} y,\label{C1}\tag{C1}\\
x\mathrel{\mathsf{C}} y\rarrow y\mathrel{\mathsf{C}} x,\label{C2}\tag{C2}\\
x\leq y\wedge z\mathrel{\mathsf{C}} x\rarrow z\mathrel{\mathsf{C}} y\,, \label{C3}\tag{C3}\\
x \con y+ z\rarrow x \con y\vee x\con z\,.\label{C4}\tag{C4}
\end{gather}
The complement of $\con$ is denoted by `$\separ$', and in the case $x\separ y$ we say that $x$ is \emph{separated from} $y$.  The class of all Boolean Contact Algebras is denoted by `$\BCA$'. The class of all \emph{complete} BCAs (i.e., these whose underlying BA is complete) is denoted by `$\BCAc$'.
Generally, for any class $\boldsymbol{\mathsf{K}}$ of algebras, $\boldsymbol{\mathsf{K}^\mathsf{c}}$ is its subclass composed of complete algebras from $\boldsymbol{\mathsf{K}}$. Given $\boldsymbol{\mathsf{K}}$ and extra constraints $(\varphi_1),\ldots,(\varphi_n)$ put upon elements of $\boldsymbol{\mathsf{K}}$, the following:
\[
\boldsymbol{\mathsf{K}}+(\varphi_1)+\ldots+(\varphi_n)
\]
denotes the subclass of $\boldsymbol{\mathsf{K}}$ compose of all structures that satisfy every $(\varphi_i)$ for $1\leqslant i\leqslant n$.

% Elements of the following, wider class of algebras:
% \[
% \BWCA\defeq\BA+\eqref{C0}+\eqref{C1}+\eqref{C2}+\eqref{C3}
% \]
% will be called \emph{Boolean weak-contact algebras}.

\newsavebox{\putwoupop}
\savebox{\putwoupop}{$\exists\mathord{\twoupop}$}
\newsavebox{\putwodownop}
\savebox{\putwodownop}{$\exists\mathord{\twodownop}$}  
\newsavebox{\deftwoupop}
\savebox{\deftwoupop}{$\mathrm{df}\,\mathord{\twoupop}$}

Unlike, e.g., \citet{Duntsch-et-al-RTBCA}, we do not assume extensionality axiom for contact:
\begin{equation}\tag{C5}\label{C5}
    (x\con a\rarrow x\con b)\rarrow a\leq b\,.
\end{equation}
nor any of its equivalent forms. 
However, later on we will consider two constraints that are equivalent to \eqref{C5}.\footnote{See \eqref{eq:E-ll} and \eqref{eq:E-gg} on page~\pageref{eq:E-ll}.}

As we see, axiom \eqref{C4} says that the contact relation distributes over the binary join operation. In the sequel, let us focus on those complete BCAs in which the contact \emph{completely} distributes over join, i.e., those that satisfy the following second-order constraint for all families $J\subseteq B$:
\begin{equation}\label{C4c}\tag{\ref{C4}$^{\mathrm{c}}$}
\textstyle  x\con\bigvee J\rarrow(\exists{y\in J})\,x\con y\,.
\end{equation}

For any region $x$ let $\con(x)$ be the set of all regions that are in contact with $x$:
\[
\con(x)\defeq\{y\in B\mid y\con x\}\,.
\]
It is easy to see that the axioms \eqref{C0}, \eqref{C3} and \eqref{C4} say that:
\begin{equation*}
    \text{$\con(x)$ is a grill.}
\end{equation*}
Further, it can be observed that \eqref{C4c} says that every $\con(x)$ is a c.p. grill:
\[
\textstyle \bigvee J\in\con(x)\rarrow(\exists y\in J)\,y\in\con(x)\,.
\]
Thus, $\con$ in the above form can be seen as a quasi-modal operator $\con\colon B\to\Gr(\frB)^\circ$, since by \eqref{C0} and \eqref{C4} respectively, we have that:
\begin{gather}
\con(\zero)=\emptyset\,,\tag{{\ref{C0}$_\star$}}\\
\con(x+y)=\con(x)\cup\con(y)\,\tag{{\ref{C4}$_\star$}}\text{ for all }x,y\in B.
\end{gather}
Additionally, it is a~completely additive operator if \eqref{C4c} holds for all $J\subseteq B$:
\begin{equation}
    \con\left(\bigvee J\right)=\bigvee_{x\in J}\con(x)\,.\tag{{\ref{C4c}$_\star$}}
\end{equation}
By Proposition~\ref{prop:cpGrills-principal-ideals} we have: 
\begin{corollary}\label{cor:co-principal}
If $\langle\frB,\mathord{\con}\rangle\in\BCAc+\eqref{C4c}$, then  $\mathord{\con}\colon B\to\Gr(\frB)^\circ$ is a~co-principal quasi-modal operator.
\end{corollary}

\subsection{Subordination and non-tangential part relations}
Every contact algebra can be expanded with the~standard binary relation of \emph{non-tangential} part (inclusion):
\begin{equation*}%\tag{$\dftt{\ll}$}
    x\ll y\iffdef x\separ -y\,.
\end{equation*}
Intuitively, $x$ is a non-tangential part of $y$ if and only if $x$ does not touch the complement of~$y$ (or, in other words, $x$ is way below $y$ or $x$ is completely surrounded by $y$). Alternatively, we can characterize non-tangential inclusion in the following way:
\begin{equation}\label{eq:ll-via-C-and-O}
x\ll y\iff\con(x)\subseteq\Overl(y)\,.
\end{equation}
It is well known that non-tangential inclusion is a special case of the subordination relation from~\citep{Bezhanishvili-G-et-al-IERGSADVD}, i.e., it has the following properties:
\begin{gather}
 0\ll 0\qquad\text{and}\qquad\one\ll\one\,,\tag{S1}\label{eq:S1}\\
 x\ll y\wedge x\ll z\rarrow x\ll y\cdot z\,,\tag{S2}\\
 x\ll y\wedge z\ll y\rarrow x+z\ll y\,,\tag{S3}\\
 x\leq y\wedge y\ll z\wedge z\leq u\rarrow x\ll u\,.\tag{S4}\label{eq:S4}
\end{gather}
The following are not standard properties of subordination, but they follow from \eqref{eq:ll-via-C-and-O} and our axioms for the contact relation:
\begin{align}
    x\ll y&{}\rarrow x\leq y\,,\label{eq:llcIngr}\tag{S5}\\
    x\ll y&{}\rarrow -y\ll-x\,,\label{eq:ll-contraposition}\tag{S6}\\
    x\ll y\wedge y\cdot z=\zero{}&{}\rarrow x\notcon z\,,\label{eq:ll-below-separated}\tag{S7}
\end{align}
and together with the former four axioms they characterize non-tangential inclusion.

Subordination is definitionally equivalent to \emph{pre-contact}, which is the relation characterized by the axioms for the contact relation without symmetry, i.e., \eqref{C2}. Non-tangential inclusion is equivalent to contact.

After \cite{Roeper-RBT}, we call \emph{isolated} every region that is separated from its own complement, i.e., any region $x$ such that $x\ll x$.

Let us observe that by means of the non-tangential part of relation we can express \eqref{C4c} in an~alternative form for all $J\subseteq B$:
\begin{equation}\label{C4cll}\tag{\ref{C4}$^{\mathrm{c}}_{\mathord{\ll}}$}
 \textstyle   (\forall y\in J)\,x\ll y\rarrow x\ll\bigwedge J\,.
\end{equation}

\section{Examples and counterexamples}\label{sec:examples}

From now on we focus on the class of \emph{complete} Boolean Contact Algebras. For brevity let us introduce the following definition:
\begin{equation*}%\tag{$\dftt{\Cfourc}$}
\Cfourc\defeq\BCAc+\eqref{C4c}\,.
\end{equation*}
To see that this class is not empty, observe that the following is true:
\begin{proposition}\label{prop:C=O->C4c}
Let $\langle\frB,\mathord{\con}\rangle$ be a complete BCA where $\mathord{\con}=\mathord{\Overl}$, then in $\frB$ contact completely distributes over join. 
\end{proposition}
\begin{proof}
Let $J\subseteq B$. We have:
\begin{align*}
\textstyle x\Overl \bigvee J{}&\textstyle\iff x\cdot\bigvee J\neq \zero\\
&{}\iff(\exists y\in J)\,x\cdot y\neq\zero\\
&{}\iff (\exists y\in J)\,x\Overl y\,.\qedhere
\end{align*}
\end{proof}
Since the class of complete BCAs in which contact is overlap is non-empty, as it contains, e.g., the power set algebra of $\Nat$ (the set of natural numbers), the class of complete BCAs that satisfy \eqref{C4c} is non-empty as well. 

Another example of the class of contact algebras that satisfy the axiom is the class of all \emph{relational contact algebras}, initiated independently by \cite{Galton-TMODS,Galton-QSC} and \cite{Vakarelov-PML}, and later developed by \cite{Duntsch-et-al-RBTODSAPA}. Given a reflexive and symmetric frame $\langle W,\mathord{\Rel}\rangle$, i.e., a non-empty set of worlds with an accessibility relation $\mathord{\Rel}\subseteq W\times W$ that is reflexive and symmetric, we expand the power set algebra $\power(W)$ with the following relation:
\begin{equation*}
    A \conr B\iffdef (\exists x\in A)(\exists y\in B)\,x\Rel y\,.
\end{equation*}
It is routine to verify that $\conr$ satisfies axioms \eqref{C0}--\eqref{C4}, so it is a Boolean contact algebra. Moreover:
\begin{proposition}\label{prop:relational-satisfy-C4c}
Every relational contact algebra belongs to the class $\Cfourc$.
\end{proposition}
\noindent This follows from the facts that the suprema in $\power(W)$ are set-theoretic sums, and so if $A\conr\bigcup_{i\in I}B_i$, there must be $x\in A$ and $y\in B_i$ for some $i\in I$ such that $x\Rel y$.

Thanks to \cite{Duntsch-et-al-RBTODSAPA} we know that every contact algebra can be embedded in a~relational contact algebra. Because of this, algebras of this kind may be treated as the standard examples of elements of the class $\Cfourc$.

% again we see that the class of contact algebras that satisfy \eqref{C4c} is non-empty. 

% In particular, every relational Boolean contact algebra is an example of a complete atomic BCA that meets \eqref{C4c}. 

Given a topological space $X$, let $\RC(X)$ be the complete Boolean algebra of all regular closed subsets of $X$, i.e., the sets $A\in\power(X)$ for which $A=\Cl\Int A$.\footnote{$\Cl$ and $\Int$ are the standard topological operations of, respectively, closure and interior.} The Boolean operations are given by the following identities:\label{page:topological-interpretation-of-C}
\begin{align*}
    A\cdot B&{}\defeq \Cl\Int(A\cap B)\,,\\
    A+ B&{}\defeq A\cup B\,,\\
    \mathord{-}A&{}\defeq\Cl (X\setminus A)\,,\\
    \bigvee\calS&{}\defeq\Cl\bigcup\calS\,.%\\
    %\bigwedge\calS&{}\defeq\bigcap\calS\,.
\end{align*}
Details concerning Regular Closed Algebras can be found, e.g., in  \citep{Koppelberg-GTBA}.

According to a~general fact, for any topological space~$X$, $\RC(X)$ with the contact interpreted as non-emptiness of the set-theoretic intersection:
\[
A\con B\iffdef A\cap B\neq\emptyset\,,
\]
is a~Boolean Contact Algebra in which: 
\[
A\ll B\iff A\subseteq\Int B\,.
\]

Using the topological interpretation of BCAs, we can also prove a general fact that entails non-emptiness of the subclass of $\Cfourc$, whose elements satisfy $\mathord{\Overl}\subsetneq\mathord{\con}$. To this end, recall that by an \emph{Alexandroff space} we mean any topological space in which the intersection of any family of open sets is open (equivalently: the sum of any family of closed sets is closed).
\begin{proposition}
If $X$ is an Alexandroff space, then $\langle\RC(X),\mathord{\con}\rangle$ satisfies \eqref{C4c}. 
\end{proposition}
\begin{proof}
Let $X$ be an Alexandroff space. If $\calS\subseteq\power(X)$, then:
\[
\Cl\bigcup\calS\subseteq\bigcup_{S\in\calS}\Cl S\,.
\]
Therefore, in $\RC(X)$, if $A\con\bigvee\calS$, then $A\cap \bigcup_{S\in\calS}S\neq\emptyset$, which means that for some $S\in\calS$, $A\cap S\neq\emptyset$, as required.
\end{proof}

\begin{example}
We can now observe that there are Alexandroff spaces whose Regular Closed Algebras are different from the Clopen Algebras, and in which contact is not overlap. Let us consider the domain: 
\[
A\defeq(\omega^+\times\{1\})\cup(\omega^+\times\{2\})\cup\{\bot\}
\]
i.e. two standardly ordered copies of the positive integers, with $\bot\notin\omega^+$ as the bottom element (i.e., $\bot\leq x$, for all $x\in A$). Let $L\defeq(\omega^+\times\{1\})\cup\{\bot\}$ be the left branch of the poset $\langle A,\mathord{\leq}\rangle$, and let $R\defeq(\omega^+\times\{2\})\cup\{\bot\}$ be its right branch. Equip the poset with the upper Alexandroff topology. The interior of $L$ is $L\setminus\{\bot\}=\upop \langle 1,1\rangle$, and its closure is $L$, so $L$ is regular closed. Similar arguments apply to $R$. Thus, it is not hard to see that $\RC(A)=\{\zero,L,R,A\}$, $L\con R$, but $L\cdot R=\zero$ and $\CO(A)=\{\zero,A\}$.

We may, of course, repeat the construction for any finite or infinite number of copies of $\omega^+$ (or any other ordered---but nor necessarily well-ordered---set). In the finite case, we have $n$ many branches, $B_i$, for every $i\leqslant n$, and $2^n$ regular closed sets. In the infinite case, when we take some infinite cardinal $\kappa$, we have:
\[
A\defeq\bigcup_{\alpha<\kappa}(\omega^+\times\{\alpha\})\cup\{\bot\}\,.
\]
For every $\alpha<\kappa$, $B_\alpha=(\omega^+\times\{\alpha\})\cup\{\bot\}$ is an atom of $\RC(A)$, and every element in $\RC(A)$ is the set-theoretic sum of $B_\alpha$'s, and vice versa. Thus $|\RC(A)|=2^\kappa$.
\end{example}

It is clear that \eqref{C4c} entails \eqref{C4}, yet the converse implication is not true. 
We are going to consider some examples, including atomless (i.e. without atoms) and atomic Boolean algebras (i.e. such in which every region is the supremum of a~set of atoms).

\begin{example}\label{ex:RC(R)-not-C4c}
To show that \eqref{C4} does not entail \eqref{C4c} consider the set of reals, the (atomless) algebra $\langle\RC(\Real),\mathord{\con}\rangle$ and the following family of regular closed subsets of the reals:
  \[
  \calS\defeq\left\{\left[-x,x\right]\,\middle\vert\,0<x<\sqrt{2}\right\}\,.
  \]
  Since $\bigvee\calS=\Cl\bigcup\calS=\left[-\sqrt{2},\sqrt{2}\right]$, we have that: 
  \begin{equation*}
  \left(-\infty,-\sqrt{2}\right]\con\bigvee\calS\qquad\text{but\quad for all\ } S\in\calS, \left(-\infty,-\sqrt{2}\right]\notcon S\,.
  \end{equation*}
\end{example}

\begin{example}
Now, we will show that there exists a Boolean space $X$ whose algebra $\RC(X)$ with the standard topological contact is atomless and fails to satisfy \eqref{C4c}. To see this take the Cantor space $2^\omega$, i.e., the countable product of discrete space $\{0,1\}$ with the standard product topology. For every $f\in 2^\omega$ the set:
\[
B_n(f)\defeq\{g\in 2^\omega\mid(\forall k\leqslant n)\,g(k)=f(k)\}
\]
is a~local basis at~$f$. Let $\mathrm{Odd}$ be the set of all sequences that begin with an odd number of zeroes followed by $1$, and $\mathrm{Even}$ be the set of all sequences that begin with an even number of zeroes followed by $1$. Let $(0)$ be the constant zero sequence. Both $O\defeq\mathrm{Odd}\cup\{(0)\}$ and $E\defeq\mathrm{Even}\cup\{(0)\}$ are regular closed, while both  $\mathrm{Odd}$ and $\mathrm{Even}$ are regular open. We see that $O\con E$. Yet the space is zero-dimensional, so there is a family of clopen sets $\{A_i\mid i\in I\}$ such that $E=\bigvee_{i\in I}A_i$, and none of $A_i$ contains $(0)$. To see this, think about all clopen sets inside $\mathrm{Even}$: none of them has the constant zero sequence among its inhabitants, yet the supremum of these sets in $\RC(X)$ is equal to the closure of their set-theoretic sum, i.e., to $E$. Thus $O$ is in contact with the supremum of all $A_i$'s, but with none of the $A_i$'s themselves.

The example is interesting also for the reason that $\RO(2^\omega)$ has a dense subalgebra, i.e. $\CO(2^\omega)$, the algebra of all clopen subsets of the Cantor space, whose elements satisfy \eqref{C4c}. If $A$ is a clopen subset of $2^\omega$, and $\{B_i\mid i\in I\}$ is a family of such subsets, then if $A\con\bigvee_{i\in I}B_i$, then $\Cl A\cap\Cl\Int\bigcup_{i\in I}B_i\neq\emptyset$. Thus $A\cap\Cl\bigcup_{i\in I}B_i\neq\emptyset$, and in consequence $A\cap B_{i_0}\neq\emptyset$, for some $i_0\in I$, as required.\qed
\end{example}

Before we go on to the next example let us recall the following construction of extensions of the contact relation. Let $\ide$ be an ideal of a~Boolean algebra $\frB$ and let $\con$ be a contact relation on $\frB$, then we may define the \emph{ideal extension}\footnote{The notion is a version of the similar extension of contact via the so-called local contact algebra, as defined in \citep{Vakarelov-et-al-PATSRBTS} and \citep{Dimov-DVTDTCLCSCM-I}.} of $\con$ for all $x,y\in  B$ via:
\[
x\con_\ide y\iffdef x\con y\vee x,y\notin\ide\,.
\]
It is obvious that $\neg\,\zero\con_\ide x$, for all $x$. \eqref{C1} and \eqref{C2} are immediate. For \eqref{C3}, let $x,y,z\in B$. Observe that if $x\leq y$ and $z,x\notin\ide$, then, since every ideal is downward closed, $y\notin\ide$, and so $z\con_\ide y$. Finally, if $x\con_\ide y+ z$ and $x,y+ z\notin\ide$, then either $y\notin \ide$ or $z\notin \ide$. Thus $x\con_\ide y$ or $x\con_\ide z$, and \eqref{C4} holds either. 
\begin{example}\label{ex:N-infinite-in-contact}
Consider the power set algebra for the set of natural numbers, and the contact algebra $\langle\power(\Nat),\mathord{\Overl}\rangle$. The family $\Fin(\Nat)\eqdef\ide$ of all finite subsets of $\Nat$ is an ideal in the algebra, and so:
\[
M\con_\ide N\iffdef M\Overl N\vee M,N\notin\Fin(\Nat)
\]
is an ideal extension of $\Overl$. To show that \eqref{C4c} fails take the set $E$ of all even numbers and enumerate the set of odd numbers: $n_0,n_1,n_2,\ldots$\,. Obviously, $E\con_\ide\bigvee_{i\in\Nat}\,\{n_i\}$, since both sets are infinite, yet for no index~$i$, $E$ is in contact with $\{n_i\}$.\footnote{This BCA comes from \citep{Gruszczynski-NTP}, where it was presented to illustrate some properties of Grzegorczyk points from \citep{Grzegorczyk-AGWP}. Since the contact algebra from the example is a~Grzegorczyk contact algebra,  it demonstrates as well that \eqref{C4c} is independent from the second-order axioms of Grzegorczyk's for the contact relation.} This example shows that there are complete atomic algebras in which \eqref{C4c} fails.
\end{example}

\section{The axiom in extensional algebras}\label{sec:extensional}

We show now that adding the extensionality axiom for contact collapses the contact relation to the relation of overlapping. Therefore, in order to obtain interesting results for the class $\Cfourc$ we have to skip extensionality.

Standardly, what is called \emph{the extensionality axiom} for BCAs is \eqref{C5}. The three following constraints are its equivalents (in the class $\BCA$). Let $\langle\frB,\con\rangle$ be a BCA: 
  \begin{align}
  \tag{C5a}\label{C5a}
    &(\forall x\neq\one)(\exists y\neq\zero)\,x\notcon y\,,\\  
  \tag{\usebox{\putwodownop}}\label{eq:E-ll}
  x\neq\zero&{}\rarrow(\exists{y\in B})\,(y\neq\zero\wedge y\ll x)\,,\\
  \label{eq:E-gg}\tag{\usebox{\putwoupop}}
  x\neq\one&{}\rarrow(\exists{y\in B})\,(y\neq\one\wedge x\ll y)\,.
  \end{align}
% Indeed, if $x\neq\one$, then $-x\neq\zero$ and by \eqref{eq:E-ll} there exists $y\neq\zero$ such that $y\ll-x$. Thus, by \eqref{eq:ll-contraposition}, $x\ll-y$, and $-y\neq\zero$, as required.

% As mentioned earlier, the two conditions are equivalent to extensionality axiom~\eqref{C5}, and we leave it to the reader to convince themselves that it is the case.

Put:
\begin{equation*}
    \twoupop x\defeq\{y\in B\mid x\ll y\}\qquad\text{and}\qquad \twodownop x\defeq\{y\in B\mid y\ll x\}\,.
\end{equation*}
By \eqref{eq:llcIngr}, $x$ is an upper bound of $\twodownop x$. Assume that $z$ is such an upper bound, but $x\Overl-z$. So, $x\cdot -z\neq \zero$. By \eqref{eq:E-ll}, there exists $u\neq \zero$ such that $u\ll x$ and $u\leq-z$, i.e. $u\nleq z$, a~contradiction. In consequence we have:
\begin{proposition}\label{lem:sup-of-ll} %If $\frB\in\BPCA+\eqref{eq:E-ll}$ and $\twodownop x\neq\{\zero\}$, then $x=\bigvee\twodownop x$ and $x=\bigwedge\twoupop x$. Therefore,
If $\langle\frB,\con\rangle\in\BCA+\eqref{eq:E-ll}$, then for every region~$x$\/\textup{:} $x=\bigvee\twodownop x$ and $x=\bigwedge\twoupop x$.
\end{proposition}
In light of this, if both $\eqref{eq:E-ll}$ and $\eqref{C4cll}$ (the non-tangential version of \eqref{C4c}) hold in a~Boolean contact algebra $\langle\frB,\con\rangle$, then every region  in $B$  must be isolated:
\begin{corollary}\label{cor:C4s->x<<x}
  If $\langle\frB,\con\rangle\in\BCA+\eqref{eq:E-ll}+\eqref{C4c}$, then for every region $x$, $x\ll x$.
\end{corollary}

According to a standard result in the theory of Boolean Contact Algebras, every region is isolated iff overlap and contact coincide:
\begin{equation}
(\forall x\in B)\,x\ll x\quad\text{if and only if}\quad\mathord{\con}=\mathord{\Overl}\,.\label{eq:C=O}
\end{equation}
% Indeed, in the situation when every region is isolated, if $x\ext y$, then $x\leq-y$. But, by the assumption, $y$ is separated from its own complement, so $x$ must be separated from~$y$, by \eqref{C3}. The other inclusion is true in general, since overlap is the smallest contact relation.
% For the other direction, in the case when contact is overlap, separation is just incompatibility. Since every region is incompatible with its own complement, every region must separated from its own complement either.
From this and from Proposition~\ref{prop:C=O->C4c}, we obtain that if every region is isolated, the contact relation completely distributes over join. 

In consequence we obtain the following conclusion:
\begin{theorem}\label{th:C4c-isolated}
Let $\langle\frB,\con\rangle\in\BCA+\eqref{eq:E-ll}$. $\langle\frB,\con\rangle$ satisfies \eqref{C4c} iff every region of $B$ is isolated iff contact and overlap relations coincide.
\end{theorem}

\subsection{Independence}

To prove that \eqref{C4c} entails that every region is isolated we assumed \eqref{eq:E-ll}. To show that the assumption is relevant we must produce a~BCA in which the contact relation completely distributes over join, yet there are regions that are not isolated (or, equivalently, the overlap relation is properly included in the contact relation). As it turns out, finding such a~BCA is relatively easy. To see this take any Boolean algebra with at least 4 elements and extend it with the largest contact relation:
\[
x\con_{\mathrm{L}} y\iffdef x\neq\zero\wedge y\neq\zero\,.
\]
The Condition \eqref{eq:E-ll} fails, since the algebra has at least four elements, and thus there is a~region~$x\neq\zero$ such that $-x\neq\zero$, and thus every non-zero tangential part of $x$ is in contact with $-x$. If the algebra is either finite or infinite, it is equally easy to see that if $x$ is in contact with $\bigvee J$ for some family of regions $J$, then there must be a non-zero $y\in J$. And thus $x\con y$. Therefore we may conclude that:
\begin{theorem}
The sentence `$\mathord{\con}=\mathord{\Overl}$' (and so the sentence `every region is isolated') is independent from the axioms \eqref{C0}--\eqref{C3}, \eqref{C4c}.
\end{theorem}

\section{Between contact and modality}\label{sec:the-operator}

In this section, we are going to explore a close connection between contact relations and modalities, in our context. We have recalled that from a modal operator, we can define a quasi-modal operator, and therefore we can define a subordination relation, see for example \citep{Celani-QMA}. The other way around is not always true. However, we prove that a possibility operator can be defined in the class of complete Boolean Contact Algebras that satisfy~\eqref{C4c}.

\subsection{From contact to modality}
Let us begin with the following:
\begin{lemma}\label{lem:existence-of-m}
Let $\langle\frB,\con\rangle\in\BCAc$, then $\langle\frB,\con\rangle$ satisfies \eqref{C4c} iff for every region~$x$ there exists a~unique region $y$ such that $\con(x)=\Overl(y)$ and $x\ll y$. In particular, every finite Boolean Contact Algebra satisfies the latter property.
\end{lemma}
\begin{proof}
Let $x\in B$. In presence of \eqref{C4c}, $\con(x)$ is a completely prime grill, so by Proposition~\ref{prop:cpGrills-principal-ideals} there exists a~unique $y$ such that $\con(x)=\Overl(y)$. Since $-y$ is disjoint from $y$, it must be the case that $-y$ is separated from~$x$. From $\Overl$-extensionality we obtain the uniqueness of such region. The converse implication is immediate.
\end{proof}
\begin{example} 
Of course, existence of the region from Lemma~\ref{lem:existence-of-m} requires \eqref{C4c}. To see this, take $\RC(\Real)$ with the standard topological contact. In Example~\ref{ex:RC(R)-not-C4c} we have shown that $\RC(\Real)$ is a~BCA in which \eqref{C4c} fails. Fix an interval $I\defeq [u,w]$ with $u<w$.  If we take any element in which $I$ is non-tangentially included, say $A$, then this element has a~part $[a,b]$ such $[u,w]\subseteq(a,b)$. But then there is $D\ll [w,b]$, so $D$ is separated from $I$. Yet $D$ overlaps~$A$.\qed
\end{example}

Lemma~\ref{lem:existence-of-m} brings to light a~different axiomatization of the class $\Cfourc$ via axioms for complete BAs plus \eqref{C0}, \eqref{C1}, \eqref{C2} and:
\begin{equation}\label{eq:C4c-equivalent}\tag{\ref{C4}$_\circ^c$}
    (\forall x\in B)(\exists y\in B)\,\con(x)=\Overl(y)\,,
\end{equation}
in particular we have:
\begin{corollary}
  $\Cfourc=\BCAc+\eqref{eq:C4c-equivalent}$.
\end{corollary}

The uniqueness property from the Lemma~\ref{lem:existence-of-m} entails the existence of an operation $m\colon B\fun B$ such that:
\begin{equation}\tag{$\dftt{m}$}\label{df:m}
m(x)\defeq(\iota y)\,\con(x)=\Overl(y)\,.
\end{equation}
%In light of \eqref{df:Diamond} and the fact that in every $\frB$ from the class $\Cfourc$ the set $\con(x)$ is a~c.p. grill we obtain that $\con$, as an operator, is co-principal, and so:
By Corollary~\ref{cor:co-principal}, in every $\langle \frB,\con\rangle\in\Cfourc$, $\con$---as a~quasi-modal operator---is co-principal, so due to  \eqref{df:Diamond} we have that:
\begin{corollary}
Let $\langle \frB,\con\rangle\in\Cfourc$ and let $m\colon B\to B$ be the operator defined by \ref{df:m}. Then, $m$ is a modal possibility operator.

More specifically, we have that:
\begin{equation}
    m(x)=\bigwedge\twoupop x\,.
\end{equation}
\end{corollary}

\begin{proof}
Let $x\in B$. Since $\con(x)$ is a c.p. grill, $\con(x)^\complement$ is a~principal ideal, and by the same reasoning as in the proof of Proposition~\ref{prop:cpGrills-principal-ideals} we have that $\con(x)^\complement=\downop\bigvee\left(\con(x)^\complement\right)$. But then:
\[
\con(x)=\left(\downop\bigvee\{y\mid y\ll -x\}\right)^\complement=\Overl\left(-\bigvee\{y\mid y\ll -x\}\right)=\Overl\left(\bigwedge\{y\mid x\ll y\}\right)\,.
\]
\end{proof}
Obviously, we have that:
\begin{equation}\label{eq:C-by-m}
    x\con y\iff m(x)\cdot y\neq\zero\,,
\end{equation}
and so:
\begin{equation}
    x\ll y\iff m(x)\leq y\,.
\end{equation}

\begin{lemma}\label{lem:properties-of-m}
If $\langle\frB,\con\rangle\in\Cfourc$, then $m\colon B\fun B$ defined by \ref{df:m} is a completely additive modal possibility operator such that\/\textup{:} 
\begin{enumerate}
    \item $x\ll m(x)$, and so $x\leq m(x)$, i.e., $x+m(x)=m(x)$,
    \item $m(x)\leq-y\iff m(y)\leq -x$, which is equivalent to $m(-m(-x))\leq x$.
\end{enumerate}
\end{lemma}
\begin{proof}
Since $\con(\zero)=\Overl(\zero)=\emptyset$, we have that $m(\zero)=\zero$. 

For complete additivity, assume that for $J\subseteq B$, $m\left(\bigvee J\right)= y$. We have that:
\[
\Overl(y)=\con\left(\bigvee J\right)=\bigcup_{x\in J}\con(x)=\bigcup_{x\in J}\Overl(m(x))=\Overl\left(\bigvee_{x\in J}m(x)\right)\,,
\]
thus $y=\bigvee_{x\in J}m(x)$, by $\Overl$-extensionality.

\smallskip

Ad.1  Directly from Lemma~\ref{lem:existence-of-m} and properties of non-tangential inclusion.

\smallskip

Ad.2 Suppose $m(x)\leq -y$ and $m(y)\nleq-x$. From the first assumption we have that $m(x)\cdot y=\zero$, so $y\notin\Overl(m(x))=\con(x)$, i.e. $y\notcon x$. But from the second one we get that $x\in\Overl(m(y))=\con(y)$, i.e., $x\con y$, a contradiction. The second direction is analogous.

Now, since $-m(y)\leq-m(y)$, from the equivalence we get that $y\leq-m(-m(y))$, so $m(-m(y))\leq -y$. Substituting $-x$ for $y$ we get that $m(-m(-x))\leq x$, as required.

The other way round, if $x\leq -m(y)$, then $m(x)\leq m(-m(y))\leq -y$, by monotonicity of the $m$ operator.
\end{proof}

% \begin{corollary}
% If $\langle B,\mathord{\con}\rangle\in\Cfourc$ and $m$ is its modal operator from Lemma~\ref{lem:properties-of-m}, then $\mathord{\con}=\mathord{\con_m}$.
% \end{corollary}
% \begin{proof}
% Immediate by \eqref{eq:C-by-m}.
% \end{proof}

It is not hard to see that:
\begin{proposition}
Let $\langle\frB,\con\rangle\in\Cfourc$ and let $m\colon B\fun B$ be the operator defined by \ref{df:m}. Then, fixed points of $m$ are exactly the isolated regions.
\end{proposition}
\noindent so, in light of the following:
\begin{lemma}[{\citealt{Duntsch-et-al-RTBCA}}]
If $\langle\frB,\con\rangle$ is a BCA, then the set of all isolated elements, i.e. $\{x\in B\mid x\ll x\}$, in which contact is reduced to overlap is a~subalgebra of $\langle \frB,\con\rangle$. 
\end{lemma}
we have:
\begin{corollary}
  Given $\langle\frB,\con\rangle\in\Cfourc$, $A\defeq\{x\in B\mid m(x)=x\}$ in which contact is reduced to overlap is a~subalgebra of $\langle\frB,\con\rangle$.
\end{corollary}

The proposition below will be useful in the next example and in the sequel:
\begin{proposition}\label{prop:C-from-set-of-grills}
Let $\frB$ be a BA, and let $\mathbf{G}$ be a~family of its grills, then\/\textup{:}
\[
x\con_{\bfG} y\iffdef x\Overl y\vee (\exists\scrG\in\bfG)\,\{x,y\}\subseteq\scrG 
\]
is a~contact relation.
\end{proposition}
\begin{proof}
\eqref{C0}--\eqref{C2} are immediate. Further, let $x,y\in B$. If $x\con_{\bfG} y$ and $y\leq z$, then in the case $\{x,y\}\subseteq\scrG$, also $\{x,z\}\subseteq\scrG$, since every grill is upward closed. Thus $\con_{\bfG}$ satisfies \eqref{C3}. To prove \eqref{C4}, let $x,y,z\in B$. Observe that in the case $\{x,y+z\}\subseteq\scrG$, then either $y\in \scrG$ or $z\in \scrG$, and in consequence either $x\con_{\bfG} y$ or $x\con_{\bfG} z$, as required.
\end{proof}

\begin{example}\label{ex:m-is-not-closure-operator}
We can see that $m$ satisfies all but one property of a~modal closure operator: $m(m(x))\leq m(x)$. Let us observe that the property does not hold for $m$ in general. Consider the eight-element algebra $\frB$ in Figure~\ref{fig:m-not-closure}, together with the two grills $\scrG_1\defeq(\downop b)^\complement$ and $\scrG_2\defeq(\downop c)^\complement$.

\begin{figure}
\begin{tikzpicture}
\node[point,label={[label distance=-2pt]below:{\footnotesize $\zero$}}] (z) at (0,0) {};
%\node[] at (z.south) {$\zero$};
\node[point,anchor=south,label={[label distance=-2pt]below left:{\footnotesize $a$}}] (a) at (-1,1) {};
\node[point,anchor=south,label={[label distance=-2pt]below left:{\footnotesize $b$}}] (b) at (0,1) {};
\node[point,anchor=south,label={[label distance=-2pt]below right:{\footnotesize $c$}}] (c) at (1,1) {};
\node[point,anchor=south,label={[label distance=-2pt]above left:{\footnotesize $ab$}}] (ab) at (-1,2) {};
\node[point,anchor=south,label={[label distance=-2pt]above left:{\footnotesize $ac$}}] (ac) at (0,2) {};
\node[point,anchor=south,label={[label distance=-2pt]above right:{\footnotesize $bc$}}] (bc) at (1,2) {};
\node[point,anchor=south,,label={[label distance=-2pt]above:{\footnotesize $\one$}}] (o) at (0,3) {};
\draw (z) -- (a) (z) -- (b) (z) -- (c) (a) -- (ab) (b) -- (ab) (b) -- (bc) (c) -- (bc) (a) -- (ac) (c) -- (ac) (ab) -- (o) (ac) -- (o) (bc) -- (o);
\node (e) at ($ (z) !.5! (c) $) {};
\node[draw=cyan,ellipse,minimum width=2.5cm,minimum height=1cm,rotate=45] (ell) at (e) {};
\node (e2) at ($ (z) !.5! (b)$) {};
\node[draw=orange,ellipse,minimum width=2cm,minimum height=1cm,rotate=90] (ell2) at (e2) {};
\end{tikzpicture}
\caption{In general, $m$ is not a~closure operator.}\label{fig:m-not-closure}
\end{figure}
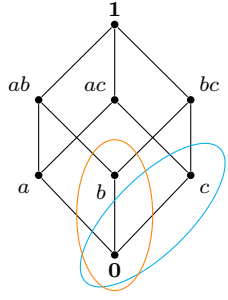

Take $\con_\bfG$ for $\bfG\defeq\{\scrG_1,\scrG_2\}$ and apply Proposition~\ref{prop:C-from-set-of-grills}. We see that:
\begin{gather*}
\con_\bfG(a)=\con_\bfG(ab)=\con_\bfG(ac)=\con_\bfG(bc)=\Overl(\one)\\
\con_\bfG(b)=\Overl(ab)\quad\text{and}\quad\con_\bfG(c)=\Overl(bc)\,.
\end{gather*}
In consequence $m(b)=ab$ and $m(m(b))=m(ab)=\one$, so $m$ in general is not a~closure operator. Observe as well that $m[B]=\{\zero,ab,bc,\one\}$ is a subboolean subalgebra\footnote{After \cite{deVries-CSC}, $\frA$ is a~\emph{subboolean subalgebra} of $\frB$ iff $A\subseteq B$ and $\frA$ is a~Boolean algebra in the partial order inherited from~$\frB$.} of $\frB$, yet it is not its~subalgebra. \qed
\end{example}

Let us observe that if $m(x)\ll m(x)$ for an element $x\in B$, then $m(x)\notcon-m(x)$, i.e.: 
\[
-m(x)\notin\con(m(x))=\Overl(m(m(x))\,.  
\]
Thus $-m(x)\leq-m(m(x))$, which means that $m(m(x))\leq m(x)$. However, Example~\ref{ex:m-is-not-closure-operator} shows that $m$ is not always a~closure operator, thus in general it is not true that every $m(x)$ is an isolated region.

The same example shows as well that in the general case $m(x)\ll m(x)$ does not entail that $x$ is isolated. For $a$ we have that $m(a)=\one$, so $m(a)\ll m(a)$, but $m(a)\neq a$, so $a$ is not isolated.

However, we have that:
\begin{proposition}
Let $\langle\frB,\con\rangle\in\Cfourc$ and let $m\colon B\fun B$ be the operator defined by \ref{df:m}. Then for all $x\in B$,
\begin{equation}
    x\notcon-x\rarrow m(x)\notcon -m(x)\,.
\end{equation}
\end{proposition}
\begin{proof}
Let $x\in B$. If $x\notcon -x$, then $-x\notin\con(x)=\Overl(m(x))$, and thus $-x\cdot m(x)=\zero$. From \eqref{eq:ll-contraposition} and Lemma~\ref{lem:properties-of-m}.1 we get that $-m(x)\notcon m(x)$.
\end{proof}

\begin{example}
The $m$ operator may not preserve Boolean operations other than those from Lemma~\ref{lem:properties-of-m}. To see this let us consider the eight-element Boolean algebra in Figure~\ref{fig:m-not-homomorphism}.
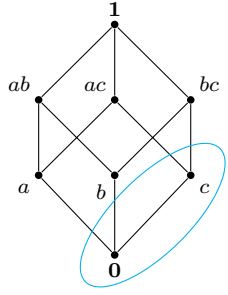
\begin{figure}
\begin{tikzpicture}
\node[point,label={[label distance=-2pt]below:{\footnotesize $\zero$}}] (z) at (0,0) {};
%\node[] at (z.south) {$\zero$};
\node[point,anchor=south,label={[label distance=-2pt]below left:{\footnotesize $a$}}] (a) at (-1,1) {};
\node[point,anchor=south,label={[label distance=-2pt]below left:{\footnotesize $b$}}] (b) at (0,1) {};
\node[point,anchor=south,label={[label distance=-2pt]below right:{\footnotesize $c$}}] (c) at (1,1) {};
\node[point,anchor=south,label={[label distance=-2pt]above left:{\footnotesize $ab$}}] (ab) at (-1,2) {};
\node[point,anchor=south,label={[label distance=-2pt]above left:{\footnotesize $ac$}}] (ac) at (0,2) {};
\node[point,anchor=south,label={[label distance=-2pt]above right:{\footnotesize $bc$}}] (bc) at (1,2) {};
\node[point,anchor=south,,label={[label distance=-2pt]above:{\footnotesize $\one$}}] (o) at (0,3) {};
\draw (z) -- (a) (z) -- (b) (z) -- (c) (a) -- (ab) (b) -- (ab) (b) -- (bc) (c) -- (bc) (a) -- (ac) (c) -- (ac) (ab) -- (o) (ac) -- (o) (bc) -- (o);
\node (e) at ($ (z) !.5! (c) $) {};
\node[draw=cyan,ellipse,minimum width=2.5cm,minimum height=1cm,rotate=45] (ell) at (e) {};
\end{tikzpicture}\caption{In general, $m$ does not have to preserve all Boolean operations.}\label{fig:m-not-homomorphism}
\end{figure}
Using the co-principal grill $\scrG\defeq(\downop c)^\complement$ we can define the contact relation:
\[
x\con_{\scrG} y\iffdef x\Overl y\vee x,y\in\scrG\,.
\]
It is routine to verify that:
\begin{gather*}
    m(c)=c\\
    m(a)=m(b)=m(ab)=ab=-c\\
    m(ac)=m(bc)=\one\,.
\end{gather*}
From these we can see that:
\[
-m(a)=c\qquad\text{but}\qquad m(-a)=m(bc)=\one
\]
and:
\[
\zero=m(a\cdot b)\neq m(a)\cdot m(b)=ab\,.
\]
Thus $m$ is not a~homomorphism. However, observe that $m[B]=\{\zero,c,-c,\one\}$ is a four element subalgebra of $B$.\footnote{The Contact Algebra from this example is a~Resolution Contact Algebra from Section~\ref{sec:RCA}.}\qed
\end{example}

We can also generalize the example in the following way.

\begin{proposition}
Let $\langle\frB,\con\rangle\in\Cfourc$ and suppose that $\frB$ is an atomic Boolean algebra. If $a\in B$ is an atom and the contact relation $\con$ is defined by\/\textup{:}
\[
x\con y\iffdef x\Overl y\vee x,y\in(\downop a)^\complement
\]
then $m[B]=\{\zero,a,-a,\one\}$.
\end{proposition}
\begin{proof}
Firstly, observe that for the atom $a$, by definition of $\con$, we have $m(a)=a$. %It is due to the fact that if $a$ were in contact with some element $x$ which is disjoint from $a$, then $a$ would have to be in contact with at least one atom of $x$, by \eqref{C4c}.

Secondly, if we take any element $x\notin\{\zero,a,\one\}$, then $x\notin\downop a$, so we have two possibilities: either (a) $x\cdot a=\zero$ or (b) $a<x$. Let $y\in B$ such that $y\neq \zero$. In (a), if $y\con x$, then $y\Overl -a$, since $x\leq-a$. If $y\Overl-a$, then $y\notin\downop a$, so $y\con x$, and so $m(x)=-a$. In (b), %if we take any non-zero $y$, 
then either $y=a$ and $y\con x$, or $y\Overl-a$, and so $y\notin\downop a$, which means that $y\con x$. So $\con(x)=\Overl(\one)$ and $m(x)=\one$.

In consequence $m[B]=\{\zero,a,-a,\one\}$.
\end{proof}

% \begin{proposition}\marginpar{\tiny Delete?}
% What if we have a BCA such that $m[B]=\{\zero,x,-x,\one\}$, for some $x\in B$?
% \end{proposition}

\begin{example}
Let us have a look at one more example. We consider the same BA reduct as in the previous two examples, but this time with contact determined by the grill $(\downop bc)^\complement$, i.e., the relation defined by:
\[
x\con y\iffdef x\Overl y\vee x,y\in(\downop bc)^\complement\,.
\]

%label={[label distance=1cm]above:\rotatebox{-45}{queen}}
\begin{figure}
\begin{tikzpicture}
\node[point,label={[label distance=-2pt]below:{\footnotesize $\zero$}}] (z) at (0,0) {};
%\node[] at (z.south) {$\zero$};
\node[point,anchor=south,label={[label distance=-2pt]below left:{\footnotesize $a$}}] (a) at (-1,1) {};
\node[point,anchor=south,label={[label distance=-2pt]below left:{\footnotesize $b$}}] (b) at (0,1) {};
\node[point,anchor=south,label={[label distance=-2pt]below right:{\footnotesize $c$}}] (c) at (1,1) {};
\node[point,anchor=south,label={[label distance=-2pt]above left:{\footnotesize $ab$}}] (ab) at (-1,2) {};
\node[point,anchor=south,label={[label distance=-2pt]above left:{\footnotesize $ac$}}] (ac) at (0,2) {};
\node[point,anchor=south,label={[label distance=-2pt]above right:{\footnotesize $bc$}}] (bc) at (1,2) {};
\node[point,anchor=south,,label={[label distance=-2pt]above:{\footnotesize $\one$}}] (o) at (0,3) {};
\draw (z) -- (a) (z) -- (b) (z) -- (c) (a) -- (ab) (b) -- (ab) (b) -- (bc) (c) -- (bc) (a) -- (ac) (c) -- (ac) (ab) -- (o) (ac) -- (o) (bc) -- (o);
\node (0c) at ($ (z) !.5! (c) $) {};
\node (bbc) at ($ (b) !.5! (bc) $) {};
\node (e) at ( $(0c) !.5! (bbc)$ ) {};
\node[draw=cyan,ellipse,minimum width=3.5cm,minimum height=2cm,rotate=45] (ell) at (e) {};
\end{tikzpicture}\caption{An algebra in which all objects are fixed points of the $m$ operator.}\label{fig:fixed-points}
\end{figure}
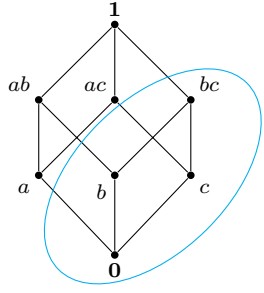
As it can be seen in Figure~\ref{fig:fixed-points}, every element is separated from its complement, thus all objects are fixed points of the operator~$m$.
\qed
\end{example}

The reason for all elements to be isolated in the example above is that if we want to have a non-empty relation $\mathord{\con}\setminus\mathord{\Overl}$ for finite structures, we must require at least two distinct atoms to be in contact. This does not have to be the case in infinite BCAs (i.e., algebras satisfying the standard \eqref{C4} axiom) as can be seen from Example~\ref{ex:N-infinite-in-contact}. However we have the following result:
\begin{proposition}
Let $\langle\frB,\con\rangle\in\Cfourc$. If $\frB$ is an infinite, atomic Boolean algebra such that $a\notcon b$ for all distinct atoms $a$ and $b$, then every element $x\in B$ is isolated, and so $\con=\Overl$.
\end{proposition}
\begin{proof}
Take $x\notin\{\zero,\one\}$ and consider the sets $A_x$ and $A_{-x}$ of all atoms that are, respectively, below $x$ and $-x$. If $x\con-x$, that is $\bigvee A_x\con\bigvee A_{-x}$, then applying \eqref{C4c} twice we obtain that there must be atoms $a\leq x$ and $b\leq-x$ such that $a\con b$.
\end{proof}

If $c$ is a co-atom of an atomic algebra, then there can only be one atom $a\notin\downop c$, so any two different atoms are isolated. In consequence, under the conditions stated above, if $c$ is a~co-atom of $\frB$ and the contact is given by\/\textup{:}
\[
x\con y\iffdef x\Overl y\vee x,y\in(\downop c)^\complement\,,
\]
then every element of $\frB$ is isolated.

% \begin{proposition}\marginpar{\tiny Problem.}
% If $x$ is an element of an atomic $\frB\in\Cfourc$ that omits $n$ atoms, then what can be said about the relation $\con$ as defined above and the operation $m$? How many elements does $m[B]$ have? 
% \end{proposition}

% \begin{proposition}\marginpar{\tiny Problem.}
% If $x$ is an element of an atomic $\frB\in\Cfourc$ with countably many atoms that omits $\aleph_0$ atoms, then what can be said about the relation $\con$ as defined above? 
% \end{proposition}

\subsection{From modality to contact}

Let us begin with the following:
\begin{definition}
Any modal algebra $\langle B,\Diamond\rangle$ whose possibility operator satisfies the following two conditions:
\begin{gather*}
    x\leq\Diamond x\,,\tag{$\mathrm{T}_{\diamond}$}\label{eq:T-Diamond}\\
    \Diamond\Box x\leq x\,,\tag{$\mathrm{B}_{\diamond}$}\label{eq:B-Diamond}
\end{gather*}
where $\Box\defeq {-}\Diamond{-}$, is called a \emph{KTB-algebra}.
\end{definition}

\begin{lemma}\label{lem:from-Diamond-to-C}
If $\langle\frB,\Diamond\rangle$ is a complete KTB-algebra, then\/\textup{:}
\begin{equation*}%\tag{$\dftt{\conm}$}\label{eq:conm}
\conm\defeq\{\langle x, y\rangle\mid x\cdot \Diamond y\neq\zero\}
\end{equation*}
is a contact relation that satisfies \eqref{C4c}. Moreover, $\Diamond=m$, where $m$ is the modal operator for~$\conm$ introduced by \eqref{df:m}.
\end{lemma}
\begin{proof}
It is obvious that $\zero\nconm x$, for any region~$x$. Reflexivity and transitivity of $\conm$ follow, respectively, from \eqref{eq:T-Diamond} and \eqref{eq:B-Diamond}. Further, if $x\cdot \Diamond(y)\neq\zero$ and $x\leq z$, then $z\cdot \Diamond(y)\neq\zero$. Finally, let $x\in B$. If $\bigvee J\in\conm(x)$ for a family $J\subseteq B$, then in light of:
\[
\conm(x)=(\downop-\Diamond x)^\complement\,,
\]
we have that there is at least one $y\in J$ for which $y\in\conm(x)$. 

By the definition of $\conm$, we have that:
\[
\conm(x)=\Overl(\Diamond x)\,.
\]
On the other hand, if we define $m$ in terms of $\conm$ by means of \eqref{df:m}, then for every $x\in B$, $m(x)$ is the only element such that $\conm(x)=\Overl(m(x))$. Therefore $\Diamond x=m(x)$ for all $x\in B$, and the two operators are equal.
\end{proof}

\subsection{An isomorphism of categories}

Now, we are going to endow the class of complete Boolean Contact Algebras that satisfy axiom \eqref{C4c} with suitable morphisms to turn them into a category. We intend to show that this category is isomorphic to the category of modal algebras and standard homomorphisms. Note that in the literature, see for example \cite{Bezhanishvili-G-et-al-IERGSADVD} and \cite{Celani-QMA}, there are at least two kinds of morphisms between Boolean algebras with a subordination relation. Our approach coincides with the definition given by \cite{Celani-QMA} where similar morphisms are called $q$-morphisms.

\begin{definition}
Given two algebras $\langle\frB_1,\con_1\rangle,\langle \frB_2,\con_2\rangle\in\Cfourc$, a~mapping $h\colon B_1\to B_2$ is a \emph{$p$-morphism}%\footnote{The idea comes from \citep{Celani-QMA}, w} 
iff it is a~homomorphism such that:
\begin{gather*}
h(x)\con_2 h(y)\rarrow x\con_1 y\,,\tag{P1}\label{eq:P1}\\
h(z)\ll_2 y\rarrow(\exists x\in B_1)(z\ll_1 x\wedge h(x)\leq_2 y)\,.\tag{P2}\label{eq:P2}
\end{gather*}
Of course, \eqref{eq:P1} is equivalent to:
\[
x\ll_1 y\rarrow h(x)\ll_2 h(y)\,.
\]
\end{definition}

\begin{lemma}
If $\langle\frB_1,\con_1\rangle,\langle \frB_2,\con_2\rangle\in\Cfourc$ and $h\colon B_1\to B_2$ is a $p$-morphism, then $h(m_1(x))=m_2(h(x))$, where $m_1$ and $m_2$ are the modal operators introduced by \eqref{df:m}. 
\end{lemma}
\begin{proof}
% Given $x\in B_2$ it is enough to show that $\con_2(h(x))=\Overl(h(m_1(x)))$, and to see this is enough to prove that $\Overl(m_2(h(x))=\Overl(h(m_1(x)))$.
Given $x\in B_2$, we know that $x\ll_1 m_1(x)$, so $h(x)\ll_2 h(m_1(x))$, which means that $h(x)\notcon_2-h(m_1(x))$. Therefore $-h(m_1(x))\notin\con_2(h(x))=\Overl(m_2(h(x)))$. Thus $-h(m_1(x))\cdot m_2(h(x))=\zero$, i.e. $m_2(h(x))\leq h(m_1(x))$.

Now, we will see that $h(m_1(x))\leq m_2(h(x))$. From $h(x)\ll_2 m_2(h(x))$, we get that there exists $y\in B_1$ such that $x\ll_1 y$ and $h(y)\leq m_2(h(x))$. So, $x\notcon_1 -y$, i.e., $-y\notin \con_1(x)=\Overl(m_1(x))$. It follows that $m_1(x)\leq y$ and therefore $h(m_1(x))\leq h(y)\leq m_2(h(x))$.
\end{proof}

\begin{lemma}
If $\langle\frB_1,\Diamond_1\rangle,\langle \frB_2,\Diamond_2\rangle\in\KTBc$ and $h\colon B_1\to B_2$ is a homomorphism of modal algebras, then $h$ satisfies properties \eqref{eq:P1} and \eqref{eq:P2} for the contact relations $\con_{\diamond_1}$ and $\con_{\diamond_2}$.
\end{lemma}
\begin{proof}
Let us prove \eqref{eq:P1} by contraposition. Suppose that $x\notcon_{\Diamond_1} y$ for $x,y\in B_1$. Then, $x\leq -\Diamond_1y$. Since $h$ is a homomorphism of modal algebras, we get $h(x)\leq -\Diamond_2(h(y))$ and thus $h(x)\notcon_{\Diamond_2}h(y)$.

Now, suppose that $h(z)\ll_2 y$ for $z\in B_1$ and $y\in B_2$. Then $h(z)\cdot \Diamond_2(-y)=\zero$ and we get $h(z)\leq \Box_2 y$. Since $\Diamond_2$ is monotone and by \ref{eq:B-Diamond}, $h(\Diamond_1z)=\Diamond_2(h(z))\leq \Diamond_2\Box_2 y\leq y$. By \ref{eq:T-Diamond}, $z\ll_1 \Diamond_1z$ and the result follows.
\end{proof}

\begin{proposition}
The class $\Cfourc$ together with $p$-morphisms form a~category with the identity functions serving as the identity morphisms. 
\end{proposition}
\begin{proof}
It is obvious that every identity is a~homomorphism and satisfies \eqref{eq:P1}. For \eqref{eq:P2}, in the case $z\ll y$, it is enough to take $y$ as $x$. Applying the definition twice we show that the composition of two p-morphisms is a~p-morphism.
\end{proof}

It is clear that the class $\KTBc$ of complete KTB algebras together with the standard homomorphisms form a category. For the remainder of this section let us treat $\Cfourc$ and $\KTBc$ as categories with suitable objects and morphisms.

From the above we can see that there is a~covariant functor $F\colon\Cfourc\to\KTBc$ which sends a~complete BCA satisfying \eqref{C4c} to a~modal algebra, and such that for every $f\in\Hom_{\Cfourc}(B_1,B_2)$, $f$ is also in  $\Hom_{\KTBc}(B_1,B_2)$, i.e. $F(f)=f$. Analogously, there is a~covariant functor $G\colon\KTBc\to\Cfourc$, such that for every $h\in\Hom_{\KTBc}(B_1,B_2)$, $G(h)=h$.

Moreover, it is the case that:
\[
G\circ F=1_{\Cfourc}\qquad\text{and}\qquad F\circ G=1_{\KTBc}\,,
\]
where $1_{\Cfourc}$ and $1_{\KTBc}$ are the identity functors for the respective categories. Thus%\marginpar{\tiny Give a short proof of this ---  if necessary.}

% It is easy to see that $F$ is full and faithful. It is also dense, in the sense that every object $C$ in $\KTBc$ is isomorphic to $F(B)$, for some $B\in\Cfourc$. To see this observe that for $\langle B,\mathord{\Diamond}\rangle\in\KTBc$, $F(\langle B,\mathord{\conm}\rangle)=\langle B,\mathord{\Diamond}\rangle$, by Lemma~\ref{lem:from-Diamond-to-C}. Thus we have:
\begin{theorem}
The categories $\Cfourc$ and $\KTBc$ are isomorphic.
\end{theorem}

\begin{figure}
\begin{equation*}
    \begin{tikzcd}[row sep=huge]
    \langle B,\con\rangle \arrow[r,"e"] \arrow[dr,"F"] & \langle B,\con,m\rangle \arrow[d,"r"]\\
    &\langle B,m\rangle
    \end{tikzcd}\qquad\quad
    \begin{tikzcd}[row sep=huge]
    \langle B,\Diamond\rangle \arrow[r,"e"] \arrow[dr,"G"] & \langle B,\conm,\Diamond\rangle \arrow[d,"r"]\\
    &\langle B,\conm\rangle
    \end{tikzcd}
\end{equation*}\caption{In both cases, $F$ and $G$ are compositions of the expansion and the reduction of suitable structures, with respect to the objects of both categories.}
\end{figure}
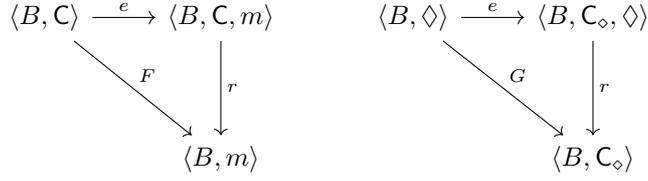

\section{Resolution contact algebras}\label{sec:RCA}

In the final section, we examine a~proper subclass of $\Cfourc$ that will serve as a~geometrical interpretation of both the contact relation that satisfies \eqref{C4c} and the modal operator defined via the contact. The inspiration for this kind of interpretation comes from \citep{Pawlak-RS}, \citep{Worboys-IIFRSD} and \citep{Duntsch-et-al-AOAR}.

\begin{definition}
A~\emph{partition} of a~Boolean algebra $\frB$ is any non-empty set $P$ of non-zero and disjoint regions of $B$ that add up to the unity: $\bigvee P=\one$.
\end{definition}

Let $\frB\in\BAc$ and let $P\defeq\{p_i\mid i\in I\}$ be a partition. We define the following relation $\conp\subseteq B\times B$:
\begin{equation*}
x\conp y\iffdef(\exists i\in I)\,(x\Overl p_i\wedge y\Overl p_i)\,.
\end{equation*}
It is routine to verify that $\conp$ is a contact relation which satisfies \eqref{C4c}. For every element $p_i$ of the partition, $\langle\downop p_i,\mathord{\con}_i\rangle$ where $\mathord{\con}_i\defeq\mathord{\conp}\cap(\downop p_i\times \downop p_i)$ is a~BCA with the full contact relation, so in particular, it satisfies \eqref{C4c}.

We adopt the following conventions: every partition of $\frB$ is called its \emph{resolution}\footnote{The name comes from~\citep{Worboys-IIFRSD}, yet unlike there we do not limit it to finite partitions.}, and the elements of the partition are called \emph{cells}. Any complete Boolean algebra expanded with $\conp$ for a~given partition $P$ is called a \emph{resolution contact algebra}. $\RCA$ is the class of such algebras. In the case $x\conp y$ we say that $x$ is in \emph{c-contact} with $y$. 

For example, the regions $x$ and $y$ in Figure~\ref{fig:coarse-p-contact} are in c-contact, since they overlap a~common cell from the sixteen-element partition. From the perspective of the picture, those regions may seem to be way apart, but we can think of the resolution as the frame of reference for the comparison of regions with respect to the $\conp$ relation. The finer the resolution, the more precise approximation of contact between regions, as we can see in Figure~\ref{fig:finer-p-contact}.

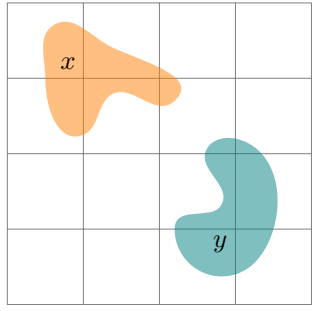
\begin{figure}
\begin{tikzpicture}
\draw[help lines] (0,0) grid (4,4);
\coordinate (a1) at (.5,3);
\coordinate (a2) at (1,2.25);
\coordinate (a3) at (1.4,2.8);
\coordinate (a4) at (2.2,2.7);
\coordinate (a5) at (1.4,3.4);
\coordinate (a6) at (0.6,3.7);
% \foreach \x in {1,2,...,6} {
% \fill (a\x) circle (1pt);
% }
\fill[orange,opacity=0.5,use Hobby shortcut] ([out angle=270]a1) .. (a2) .. ([in angle=190]a3) .. (a4) .. ([in angle=0]a5) .. (a6) .. ([in angle=90]a1) -- cycle;
\coordinate (b1) at (2.2,1);
\coordinate (b2) at (2.8,1.3);
\coordinate (b3) at (2.6,2);
\coordinate (b4) at (3,2.2);
\coordinate (b5) at (3.5,1);
\coordinate (b6) at (3,0.4);
% \foreach \x in {1,2,...,6} {
% \fill (b\x) circle (1pt);
% }
\fill[teal,opacity=0.5,use Hobby shortcut] ([out angle=90]b1) .. (b2) .. ([in angle=190]b3) .. (b4) .. ([in angle=0]b5) .. (b6) .. ([in angle=270]b1) -- cycle;
\node at (0.8,3.2) {$x$};
\node at (2.8,0.8) {$y$};
\end{tikzpicture}\caption{Regions $x$ and $y$ that are in contact with respect to\\ the sixteen-element partition.}\label{fig:coarse-p-contact}
\end{figure}

% \begin{figure}
% \begin{tikzpicture}[scale=0.25]
% \draw[help lines] (-6,-6) grid (10,10);
% \coordinate (a1) at (.5,3);
% \coordinate (a2) at (1,2.25);
% \coordinate (a3) at (1.4,2.8);
% \coordinate (a4) at (2.2,2.7);
% \coordinate (a5) at (1.4,3.4);
% \coordinate (a6) at (0.6,3.7);
% % \foreach \x in {1,2,...,6} {
% % \fill (a\x) circle (1pt);
% % }
% \fill[orange,opacity=0.5,use Hobby shortcut] ([out angle=270]a1) .. (a2) .. ([in angle=190]a3) .. (a4) .. ([in angle=0]a5) .. (a6) .. ([in angle=90]a1) -- cycle;
% \coordinate (b1) at (2.2,1);
% \coordinate (b2) at (2.8,1.3);
% \coordinate (b3) at (2.6,2);
% \coordinate (b4) at (3,2.2);
% \coordinate (b5) at (3.5,1);
% \coordinate (b6) at (3,0.4);
% % \foreach \x in {1,2,...,6} {
% % \fill (b\x) circle (1pt);
% % }
% \fill[teal,opacity=0.5,use Hobby shortcut] ([out angle=90]b1) .. (b2) .. ([in angle=190]b3) .. (b4) .. ([in angle=0]b5) .. (b6) .. ([in angle=270]b1) -- cycle;
% %\node at (0.8,3.2) {$x$};
% %\node at (2.8,0.8) {$y$};
% \end{tikzpicture}\caption{Regions $x$ and $y$ zoomed-out.}\label{fig:fine-p-contact}
% \end{figure}
Thus, the fineness of the partition is a~counterpart of the precision with which we can discern regions and their mutual relations. If we cover the space from Figure~\ref{fig:coarse-p-contact} with a finer partition, then we have a~more precise notion of contact, and more precise approximations of regions via cells.

\begin{figure}
\begin{tikzpicture}
\draw[help lines,step=0.2] (0,0) grid (4,4);
\coordinate (a1) at (.5,3);
\coordinate (a2) at (1,2.25);
\coordinate (a3) at (1.4,2.8);
\coordinate (a4) at (2.2,2.7);
\coordinate (a5) at (1.4,3.4);
\coordinate (a6) at (0.6,3.7);
% \foreach \x in {1,2,...,6} {
% \fill (a\x) circle (1pt);
% }
\fill[orange,opacity=0.5,use Hobby shortcut] ([out angle=270]a1) .. (a2) .. ([in angle=190]a3) .. (a4) .. ([in angle=0]a5) .. (a6) .. ([in angle=90]a1) -- cycle;
\coordinate (b1) at (2.2,1);
\coordinate (b2) at (2.8,1.3);
\coordinate (b3) at (2.6,2);
\coordinate (b4) at (3,2.2);
\coordinate (b5) at (3.5,1);
\coordinate (b6) at (3,0.4);
% \foreach \x in {1,2,...,6} {
% \fill (b\x) circle (1pt);
% }
\fill[teal,opacity=0.5,use Hobby shortcut] ([out angle=90]b1) .. (b2) .. ([in angle=190]b3) .. (b4) .. ([in angle=0]b5) .. (b6) .. ([in angle=270]b1) -- cycle;
\node at (0.8,3.2) {$x$};
\node at (2.8,0.8) {$y$};
\end{tikzpicture}\caption{Regions $x$ and $y$ are no longer in contact if we take a~finer partition as the frame of reference.}\label{fig:finer-p-contact}
\end{figure}
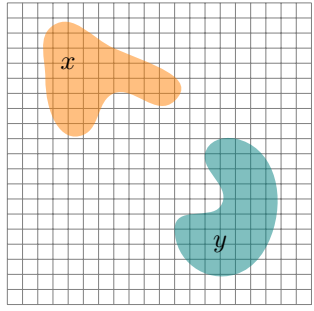

It follows from the definition of $\conp$ that every element $p_i$ of the partition $P$ indexed by $I$ must be isolated, thus $m(p_i)=p_i$. It is also the case that an arbitrary join of elements of the partition must be isolated. Indeed, let $J\subsetneq I$ and $K\defeq I\setminus J$. Consider the regions $\bigvee_{j\in J}p_j$ and $\bigvee_{k\in K}p_k$. If they are in contact, then applying \eqref{C4c} twice we get that there are $j\in J$ and $k\in K$ such that $p_j\conp p_k$, which means that they overlap a~common element from $P$. But this is only possible if $p_j=p_k$, a~contradiction. In consequence, if $J\subseteq I$, then:
\[
\textstyle m\left(\bigvee_{j\in J}p_j\right)=\bigvee_{j\in J}p_j\,.
\]
%It follows that $m[B]$ must be a~complete, atomic algebra, whose atoms are elements of the partition.

Further, for any region $x$ define $\Overl_P(x)\defeq\{p_i\in P\mid x\Overl p_i\}$, the set of all elements of the partition that overlap~$x$. For any $x$, $\Overl_P(x)$ is called \emph{the cell cover} (or \emph{c-cover}) of~$x$. 

%Analogously, the set of all elements of the partion $P$ that are parts of $x$ will be called \emph{the cell interior}. 

We have that $m(x)=\bigvee\Overl_P(x)$, and so if $x\leq p_i$, for some $i\in I$, then $m(x)=p_i$. In particular, for any $p_i$ which is not an atom, none of its proper parts is a~fixed point of the $m$ operator.

\begin{figure}
\begin{tikzpicture}
\draw[help lines] (-1,0) grid (4,5);
\coordinate (a1) at (.5,3);
\coordinate (a2) at (1,2.25);
\coordinate (a3) at (1.4,2.8);
\coordinate (a4) at (2.2,2.7);
\coordinate (a5) at (1.4,3.4);
\coordinate (a6) at (0.6,3.7);
% \foreach \x in {1,2,...,6} {
% \fill (a\x) circle (1pt);
% }
\fill[gray,opacity=0.5] (0,2) rectangle (3,4);
\fill[orange,opacity=0.5,use Hobby shortcut] ([out angle=270]a1) .. (a2) .. ([in angle=190]a3) .. (a4) .. ([in angle=0]a5) .. (a6) .. ([in angle=90]a1) -- cycle;
\node at (.8,3.2) {$x$};
\end{tikzpicture}\caption{The cell cover of the region $x$, whose cell interior\\ is empty.}\label{fig:orbit}
\end{figure}
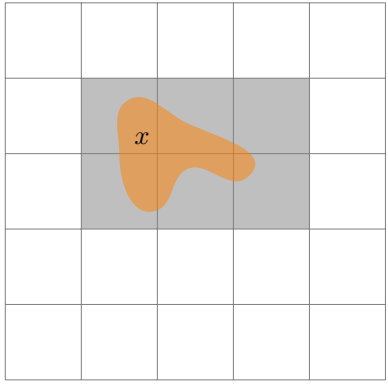

\begin{proposition}\label{prop:C_P-is-possible-contact}
Let $\langle\frB,\con_P\rangle$ be a~resolution algebra  with a partition $P$. If $x,y\in B$, then $x\con_P y$ iff the c-covers of $x$ and $y$ overlap, i.e. $m(x)\cdot m(y)\neq\zero$.
\end{proposition}
\begin{proof}
The left-to-right implication is true in general, and stems from \eqref{df:m} and the fact that $x\leq m(x)$ for every region~$x$.

For the right-to-left implication let $x,y\in B$ and assume that $m(x)\cdot m(y)\neq\zero$, so there must be a~cell $p_i$ that overlaps both $m(x)$ and $m(y)$. But this means that $p_i$ overlaps a~cell $p_j$ from the c-cover of $x$, and a~cell $p_k$ from the c-cover of~$y$. Since all three are elements of a resolution, they must be equal: $p_i=p_j=p_k$.
\end{proof}
\begin{remark}
From the above proposition we can see that $\conp$ is a~counterpart of a~\emph{possible} contact relation from \citep{Duntsch-et-al-AOAR}. However, the right-to-left implication from Proposition~\ref{prop:C_P-is-possible-contact} is not always true about the $m$ operator (in the sense that it is not true in every algebra from the class $\Cfourc$). This can be seen in a~model from Example~\ref{ex:m-is-not-closure-operator} where we have that $m(a)=\one$ and $m(c)=bc$, so $m(a)\cdot m(c)\neq\zero$. Yet $a\notcon_{\mathbf{G}}c$.
\end{remark}

Let us make a~notion of the \emph{fineness} of a~partition precise, and let us observe that the finer the partition, the fewer regions are in contact. 
\begin{definition}
Let $\frB$ be a complete Boolean algebra. Given two partitions, $P$ and $P'$, $P$ is \emph{finer} than $P'$ (in symbols: $P\unlhd P'$) iff for every $p\in P$ there is a~$p'$ in $P'$ such that $p\leq p'$. $P$ is \emph{strictly finer} than $P'$ iff $P\unlhd P'$ and there is a $p\in P$ such that $p\notin P'$.
\end{definition}

\begin{proposition}

Let $\frB$ be a complete Boolea algebra. Given a~sequence of partitions $\langle P_\alpha\mid\alpha<\kappa\rangle$ (with $\kappa$ being an~ordinal) such that $P_{\alpha+1}\unlhd P_\alpha$, let us consider the contact relation $\con_\alpha$ determined by every $P_\alpha$. We have that:
\[
\beta<\alpha\rarrow\con_{\alpha}\subseteq\con_{\beta}
\]
and if $P_\alpha$ is strictly finer than $P_\beta$, then the inclusion is proper.
\end{proposition}
\begin{proof}
Let $x,y\in B$ and suppose that $\beta<\alpha$. If $x\con_\alpha y$, then there is a $p_i\in P_\alpha$ such that $p_i\in\Overl_{P_\alpha}(x)\cap\Overl_{P_\alpha}(y)$. Since $P_\alpha$ is finer that $P_\beta$, there is $p_j\in P_\beta$ such that $p_i\leq p_j$. So $p_j\in \Overl_{P_\beta}(x)\cap\Overl_{P_\beta}(y)$, and in consequence %$m(x)\cdot m(y)\neq\zero$, so by Proposition~\ref{prop:C_P-is-possible-contact},
$x\con_\beta y$.

Suppose now that $P_\alpha$ is strictly finer than $P_\beta$, i.e. there exists a $p_i\in P_\alpha\setminus P_\beta$. Thus in $P_\beta$ there is $p_k$ such that $p_i<p_k$. In consequence there is a~non-zero region~$x$ such that $x<p_k$ and $x\cdot p_i=\zero$. Therefore $x\notcon_\alpha p_i$, but $x\con_\beta p_i$, as both $x$ and $p_i$ are parts of the same cell $p_k\in P_\beta$.
\end{proof}

Since the overlapping relation is the smallest contact relation, it must be the case that for $\langle P_\alpha\mid\alpha<\kappa\rangle$ as above:
\[
\mathord{\Overl}\subseteq\bigcap_{\alpha<\kappa}C_\alpha,
\]
and in general, it does not have to be the case that the intersection is exactly overlap. For example, look at the regular closed sets in $[0,1]^2$ considered as a~subspace of~$\Real^2$ with the standard topology. If we start with the sixteen-element partition from Figure~\ref{fig:overl-neq-intersection}, the two triangles are in c-contact in every refinement generated by quadratic subdivisions of cells. Therefore, the intersection of all contact relations $\langle \con_n\mid n<\omega\rangle$ generated by the partitions contains the overlap as a proper subset.

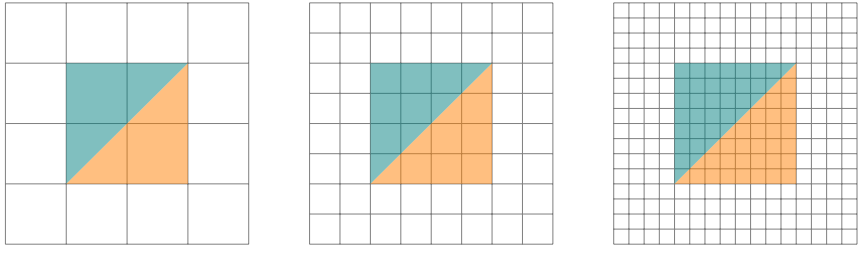
\begin{figure}
\begin{tikzpicture}[scale=0.8]
\draw[help lines,step=1] (0,0) grid (4,4);
\fill[orange,opacity=0.5] (1,1) -- (3,1) -- (3,3) -- cycle;
\fill[teal,opacity=0.5] (1,1) -- (1,3) -- (3,3) -- cycle;
\begin{scope}[xshift=5cm]
\draw[help lines,step=0.5] (0,0) grid (4,4);
\fill[orange,opacity=0.5] (1,1) -- (3,1) -- (3,3) -- cycle;
\fill[teal,opacity=0.5] (1,1) -- (1,3) -- (3,3) -- cycle;
\end{scope}
\begin{scope}[xshift=10cm]
\draw[help lines,step=0.25] (0,0) grid (4,4);
\fill[orange,opacity=0.5] (1,1) -- (3,1) -- (3,3) -- cycle;
\fill[teal,opacity=0.5] (1,1) -- (1,3) -- (3,3) -- cycle;
\end{scope}
\end{tikzpicture}\caption{The two triangles are in contact for every contact generated by any partition that is a~refinement of the initial sixteen-element partition.}\label{fig:overl-neq-intersection}
\end{figure}

However, the intersection of this kind may not be an approximation of the standard topological contact in $\RC([0,1]^2)$, which can be considered as a weak point of contact relations generated by partitions (in the case we want to treat them as approximations). For example, no two different cells from the sixteen-element partition from Figure~\ref{fig:overl-neq-intersection} will ever be in contact with respect to contact relations generated by its refinement, yet as regions of  $[0,1]^2$ (i.e., its regular closed subsets) adjacent cells have non-empty set-theoretical intersections, so they are in contact under the standard topological interpretation  of contact (see page~\pageref{page:topological-interpretation-of-C}).

\subsection{Resolution algebras and frames of ultrafilters}

\begin{proposition}
Let $\langle \frB,\con_P\rangle$ be a resolution algebra. Let $m\colon B\fun B$ be the operator defined by \eqref{df:m}. Then, $m(m(x))=m(x)$ for all $x\in B$, i.e. in every resolution algebra $m$ is a closure operator.
\end{proposition}
\begin{proof}
Let $x\in B$. We apply the characterization of $m(x)$ for resolution contact algebras plus the facts that $m$ is completely additive and that every element of the partition is a~fixed point of $m$:
\begin{align*}
    m(m(x))=m\left(\bigvee\Overl_P(x)\right)=\bigvee\{m(p_i)\mid x\Overl p_i\}=\bigvee\{p_i\mid x\Overl p_i\}=m(x)\,.
\end{align*}
\end{proof}

In light of the above, for any $\langle\frB,\con_P\rangle\in\RCA$, the standard Kripke relation $R$ on the set of ultrafilters of $\frB$ is an equivalence relation that partitions the set of ultrafilters into equivalence classes. For any ultrafilter $\ult$, let $[\ult]$ be its equivalence class and let $\Ult \frB/_R$ be the partition of $\Ult \frB$.

\begin{theorem}
If $\langle \frB,\con_P\rangle\in\RCA$ with a finite resolution $P=\{p_i\mid i\leqslant n\}$ for some $n\in\Nat$, then the Kripke relation on the set $\Ult \frB$ is an equivalence relation and there is a~one-to-one correspondence $f\colon P\to\Ult \frB/_R$ between cells and equivalence classes of ultrafilters.
\end{theorem}
\begin{proof}
Let us consider $R$ and the set $\Ult \frB/_R$. For every $p_i\in P$ define $f(p_i)\defeq s(p_i)$, the standard Stone mapping for elements of the partition. The mapping must be injective, since elements of the partition are pairwise disjoint. Further, if $s(p_i)\cap s(p_j)\neq\emptyset$, then $p_i=p_j$, and thus $s(p_i)=s(p_j)$. Still further, every ultrafilter must be in some $s(p_i)$. This follows from the fact that $\bigvee P=\one$ and $P$ is finite.

It remains to show that for every equivalence class $E$ from $\Ult \frB/_R$ there is $p_i$ such that $s(p_i)=E$. Take $\ult$ to be a~representative of $E$, and let $p_i$ be this unique element of the partition that sits inside $\ult$. This immediately entails that $\ult\in s(p_i)$. For the other direction, let $p_i\in\scrK$. If $x\in\ult$, then $x\cdot p_i\neq\zero$, i.e. $p_i\in\Overl_P(x)$. Therefore $p_i\leq\bigvee\Overl_P(x)=m(x)$, and in consequence $m(x)\in\scrK$. So $\ult\subseteq m^{-1}[\scrK]$, and so $\scrK\in[\ult]$, as required.
\end{proof}

In general, if we have a complete S5-algebra $\langle \frB,\Diamond\rangle$, its expansion $\langle \frB,\Diamond,\mathord{\con_\diamond}\rangle$ does not have to be a~resolution algebra; that is, there may be no partition~$P$ of $\frB$ such that $\mathord{\conp}= \mathord{\con_\diamond}$. For example, take as $\frB$ any atomless algebra with $\Diamond$ as a~fixed-point operator. Then $\Diamond$ is an S5 operator and  in consequence we have:
\[
x\con_\diamond y\iff x\cdot\Diamond y\neq\zero\iff x\cdot y\neq\zero\,.
\]
Yet, the algebra is atomless, so for any partition $P$ of $\frB$, any element $p_i$ of $P$ must have disjoint non-zero parts $x$ and $y$, which means that $x\notcon_\Diamond y$. Yet, as both $x$ and $y$ overlap $p_i$, it must be the case that $x\con_P y$. Therefore $\mathord{\conp}\neq\mathord{\con_\diamond}$.

But we have that:
\begin{theorem}
Given an S5 modal algebra $\langle \frB,\Diamond\rangle$, its expansion $\langle \frB,\Diamond,\mathord{\con_\diamond}\rangle$ can be embedded into a modal expansion of a~resolution algebra. 
\end{theorem}
\begin{proof}
Let $R$ be the accessibility relation on the set of all ultrafilters of $\frB$:
\[
\ult_1\Rel\ult_2\iffdef\ult_1\subseteq\Diamond^{-1}[\ult_2]\,.
\]
Since $\langle\frB,\Diamond\rangle$ is an S5-algebra, $R$ is an equivalence relation. By means of this take $\conr$ to be the contact relation between sets of ultrafilters, i.e.:
\[
A\conr D\iffdef(\exists\ult_1\in A)(\exists\ult_2\in D)\,\ult_1\Rel\ult_2\,.
\]
Since $P\defeq\Ult \frB/_R$ is a~partion of $\power(\Ult \frB)$, the relation:
\[
A\conp D\iffdef(\exists p\in P)\,A\cap p\neq\emptyset\neq D\cap p
\]
is a~contact relation that satisfies \eqref{C4c}. We have that $\mathord{\conr}=\mathord{\conp}$, since in the special case of $R$ we have:
\[
\ult_1\Rel\ult_2\iff(\exists p\in P) p=[\ult_1]=[\ult_2]\,.
\]

We now take the expansion of $\langle\power(\Ult \frB),\mathord{\conr}\rangle$ to $\frP\defeq\langle\power(\Ult \frB),m,\mathord{\conr}\rangle$, and consider the Stone mapping $s\colon B\to\power(\Ult \frB)$. In particular, we will see that:
\[
x\con_\Diamond y\iff s(x)\conr s(y)\,.
\]
Let $x,y\in B$ such that $x\con_\Diamond y$, i.e., $x\cdot \Diamond y\neq \zero$ and there exists a $\ult_1\in \Ult \frB$ that $x,\Diamond y\in \ult_1$. Then, there exists a $\ult_2\in \Ult \frB$ that $\ult_1\Rel \ult_2$ and $y\in \ult_2$. So $\ult_1\in s(x)$, $\ult_2\in s(y)$ and we get $s(x)\con_R s(y)$.
On the other hand, if $s(x)\con_R s(y)$, there exist $\ult_1,\ult_2\in \Ult \frB$ such that $\ult_1\Rel\ult_2$, $x\in\ult_1$ and $y\in\ult_2$. So, $x,\Diamond y\in \ult_1$ and we get $x\cdot \Diamond y\neq \zero$. Therefore, $x\con_\Diamond y$.
\end{proof}

\subsection*{Acknowledgements} We would like to thank to an anonymous referee for their criticism that helped to improve the paper.

% \begin{theorem}\marginpar{\tiny Hypothesis} If the canonical extension of a modal algebra~$\frB$ is a~resolution algebra, then $\frB$ is an S5-algebra. 

% \end{theorem}

% \begin{theorem}\marginpar{\tiny Hypothesis} If~$\frB=\langle B,\mathord{\Diamond}\rangle$ is a~modal algebra, and $\langle B,\mathord{\Diamond},\mathord{\con_\Diamond}\rangle$ is its expansion that can be embeeded into a~modal expansion of a~resolution algebra, then $\frB$ is an S5-algbera.

% \end{theorem}

\input{biblio.bbl}

\bibliographystyle{apalike}
\end{document}

%% file: basic_mac.tex
\newtheorem{theorem}{Theorem}[section]
\newtheorem{lemma}[theorem]{Lemma}
\newtheorem{proposition}[theorem]{Proposition}
\newtheorem{corollary}[theorem]{Corollary}

\theoremstyle{remark}
\newtheorem{definition}{Definition}[section]

\newtheorem{example}{Example}[section]
\newtheorem{remark}{Remark}[section]

\newcommand\rarrow{\longrightarrow}

\renewcommand\iff{\longleftrightarrow}
\newcommand\fun{\ensuremath{\longrightarrow}}

\newcommand\Real{\varmathbb{R}}
\newcommand\Nat{\varmathbb{N}}

\DeclareMathOperator{\Fin}{Fin}

\newcommand{\power}{\mathcal{P}} %power set
\newcommand\iffdef{\;\mathrel{\mathord{:}\mathord{\longleftrightarrow}}\;}
\newcommand\defeq{\coloneqq}
\newcommand\eqdef{\eqqcolon}

\newcommand{\zero}{\mathbf{0}} %zero of an algebra
\newcommand{\one}{\mathbf{1}} %unity of an algebra

\newcommand{\ext}{\mathrel{\bot}} %incompatibility relation

\newcommand{\fil}{\mathscr{F}} % a filter
\newcommand{\ult}{\mathscr{U}} % an ultrafilter
\newcommand{\ide}{\mathscr{I}} %an ideal

\DeclareMathOperator{\Ult}{Ult} %all ultrafilters of the algebra
\DeclareMathOperator{\Gr}{Gr} %all ideals of the algebra
\DeclareMathOperator{\Ide}{Id} %all ideals of the algebra
\DeclareMathOperator{\Hom}{Hom} %the collection of all morphisms

\DeclareMathOperator{\CO}{CO} %clopen algebra
\DeclareMathOperator{\RO}{RO} %regular open
\DeclareMathOperator{\RC}{RC} %regular closed

\newcommand\BA{\boldsymbol{\mathsf{BA}}} %the class of Boolean algebras
\newcommand\BAc{\boldsymbol{\mathsf{BA^c}}} %the class of Boolean contact algebras
\newcommand\BCA{\boldsymbol{\mathsf{BCA}}} %the class of Boolean contact algebras
\newcommand\BCAc{\boldsymbol{\mathsf{BCA^c}}} %the class of complete Boolean contact algebras
\newcommand\RCA{\boldsymbol{\mathsf{RCA}}} %the class of resolution contact algebras
\newcommand\KTBc{\boldsymbol{\mathsf{KTB^c}}} %the class of complete KTB modal algebras
\DeclareMathOperator{\con}{\mathsf{C}} %relacja stycznoci

\DeclareMathOperator{\Overl}{\mathsf{O}} %the set of all regions that overlap given object

\newcommand\mathbackslash{\raisebox{.4pt}{\texttt{/}}}
\def\notcon{% the separation relation
  \renewcommand\stacktype{L}\mathrel{\ensurestackMath{%
  \ThisStyle{\stackon[0pt]{\SavedStyle\con}{\SavedStyle\mathbackslash}}}}%
}
\let\separ=\notcon

\DeclareMathOperator{\upop}{\uparrow} %upper arrow operator
\DeclareMathOperator{\downop}{\downarrow} %lower arrow operator
\newcommand\twodownop{\raisebox{-1pt}{$\mathop{\rotatebox{90}{$\twoheadleftarrow$}}$}} %downward << closure operation
\newcommand\twoupop{\raisebox{-1pt}{$\mathop{\rotatebox{90}{$\twoheadrightarrow$}}$}} %upward << closure operation

\DeclareMathOperator{\Cl}{Cl}
\DeclareMathOperator{\Int}{Int}

\newcommand\calS{\mathcal{S}}

\newcommand\mfr{\mathfrak}
\newcommand\frA{\mfr{A}}
\newcommand\frB{\mfr{B}}

\newcommand\frP{\mfr{P}}

\newcommand\mbf{\mathbf}

\newcommand\bfG{\mbf{G}}

\newcommand{\scrG}{\mathscr{G}}

\newcommand{\scrJ}{\mathscr{J}}
\newcommand{\scrK}{\mathscr{K}}

\newcommand{\scrS}{\mathscr{S}}

\newcommand\Rel{\mathrel{R}}

\newcommand\dftt{\mathtt{df}\,}
\newcommand\suchthat{\,\middle\vert\,} %abstraction operator

%% file: biblio.bbl
\providecommand{\noop}[1]{}